\documentclass[11pt,a4paper]{amsart}
\usepackage{amsmath, amssymb}
\usepackage{amsthm} 
\usepackage[utf8]{inputenc}
\usepackage{graphicx}
\usepackage{rotating}
\usepackage[export]{adjustbox}
\graphicspath{ {./images/} }
\usepackage{color}
\usepackage{hyperref}
\usepackage{algorithm}
\usepackage{algpseudocode}
\makeatletter
\renewcommand\subsection{\@startsection{subsection}{2}%
  \z@{.7\linespacing\@plus.7\linespacing}{.5\linespacing}%
  {\normalfont\bfseries}}
\makeatother

\theoremstyle{plain}
\newtheorem{theorem}{Theorem}[section]
\newtheorem{lemma}[theorem]{Lemma}
\newtheorem{proposition}[theorem]{Proposition}

% ...

\theoremstyle{definition}
\newtheorem{definition}{Definition}[section]

\theoremstyle{remark}
\newtheorem{remark}{Remark}[section]

\title[Simulation of conformal mappings]{Efficient simulation of mixed boundary value problems and conformal mappings}

\author{Qiansheng Han}
\address{Qiansheng Han\\
Department of Mathematics with Computer Science\\
Guangdong Technion-Israel Institute of Technology\\
Shantou, Guangdong 515063\\
P.R. of China}
\email{han07024@gtiit.edu.cn}

\author{Antti Rasila\textsuperscript{*}}
\address{Antti Rasila\\
Department of Mathematics with Computer Science\\
Guangdong Technion -- Israel Institute of Technology\\
Shantou, Guangdong 515063, P.R. of China
and
Department of Mathematics\\
Technion -- Israel Institute of Technology\\
Haifa 32000, Israel}
\email{antti.rasila@gtiit.edu.cn; antti.rasila@iki.fi}
\thanks{\textsuperscript{*}Corresponding author}

\author{Tommi Sottinen}
\address{Tommi Sottinen\\
University of Vaasa\\
School of Technology and Innovations \\
P.O.Box 700\\
FIN-65101 Vaasa\\
Finland}
\email{tommi.sottinen@uwasa.fi}



\begin{document}

\begin{abstract}
In this paper, we present a stochastic method for the simulation of Laplace's equation with a mixed boundary condition in planar domains that are polygonal or bounded by circular arcs. We call this method the Reflected Walk-on-Spheres algorithm. The method combines a traditional Walk-on-Spheres algorithm with use of reflections at the Neumann boundaries. We apply our algorithm to simulate numerical conformal mappings from certain quadrilaterals to the corresponding canonical domains, and to compute their conformal moduli. Finally, we give examples of the method on three dimensional polyhedral domains, and use it to simulate the heat flow on an L-shaped insulated polyhedron.
\end{abstract}

\maketitle

\section{Introduction}

Conformal mappings are functions that preserve angles (directions and magnitudes) between two smooth curves at their point of intersection. They play an essential role in classical complex analysis, admitting a convenient characterization that they have a non-vanishing complex derivative. It is further assumed that conformal mappings are one-to-one and onto. Conformal mappings play an important role in engineering mathematics, e.g., aerospace engineering, fluid dynamics, electronics, and geodesy \cite{CalixtoAlvarenga2011,GONZALEZMATESANZ20061432,SchinzingerLaura2003,VermeerRasila2020}. They have also been recently applied in topics such as computer graphics \cite{BowersHurdal2003,KharevychSpringbornSchroder2006} and brain mapping \cite{wang2008brain}.

In this paper we consider the computation of conformal mappings \( f \) from a bounded and simply-connected domain $\Omega$ onto \( R_h = \{ z \in \mathbb{C} : 0 < \text{Re} (z) < 1, \; 0 < \text{Im} (z) < h \},\; h>0 \), which is called the canonical rectangle. The existence of such conformal mappings is guaranteed by Riemann's mapping theorem. However, this result does not provide us with an explicit formula of the mapping, which can only be derived analytically for some special cases of $\Omega$.

The problem of finding numerical representations of Riemann mappings is classical. Popular methods include the Schwarz--Christoffel method \cite{DriscollTrefethen2002}, discrete Ricci flow \cite{GuYau2008}, Wegmann's method \cite[p. 405]{Wegmann}, and circle packing \cite{Stephenson2005}. The conjugate function method used as the basis of our investigation was first introduced in \cite{HAKULA2013340}. A version of this algorithm for multiply connected domains was developed in \cite{HQR2}. It should be noted that the $hp$-FEM implementation of the conjugate function method in \cite{HAKULA2013340} works also with unbounded planar domains \cite{HRV2} and ones with strong singularities \cite{HRV3}.

We propose a new stochastic algorithm for numerical approximation of a conformal mapping of a domain onto the canonical rectangle. This approach requires solving a Dirichlet--Neumann mixed boundary value problem on the domain.  For this purpose, a variant of the Walk-on-Spheres (WoS) algorithm is used. This is a Monte Carlo method specifically designed for efficient simulation of Dirichlet type boundary value problems. Our algorithm is thus modified for solving the mixed boundary value problem. The algorithm is modified to make use of Schwarz' reflections at the Neumann boundary. We will explain this idea in detail in Section \ref{sec-refl}. To our knowledge, this is the first time attempt to develop a stochastic algorithm for numerical conformal mappings. For a survey of other currently available methods see e.g. \cite[pp. 8--11]{Kythe}. Recently, there has been substantial work on boundary integral methods for computing numerical conformal mappings \cite{Nasser1} as well as closely related problems of  conformal modulus of a quadrilateral \cite{Nasser3} and condenser capacity \cite{Nasser2}.

This paper is organized as follows. We will first present the fundamental concepts required for understanding our approach in Section \ref{sec-preli}, followed by a detailed explanations of the conjugate function method in Section \ref{sec-conj}. The theory behind the Walk-on-Sphere method is discussed in Section \ref{sec-walk}, and its variant, a reflected Walk-on-Spheres method, is introduced in Section \ref{sec-refl}. In Section \ref{sec-ex}, we give several simulation examples of numerical conformal mappings using this method, and finally discuss generalization of this approach to the three dimensional setting in Section \ref{sec-3d}.

\section{Preliminaries}\label{sec-preli}

In this section, we discuss conformal modulus of a quadrilateral, which is a concept originating from geometric function theory (see e.g. \cite{Ahlfors2006}). First recall the Riemann mapping theorem: 

\begin{theorem}[Riemann Mapping Theorem]
Let $\mathbb{D}$ be the unit disk. For a simply connected domain $\Omega$ in the plane, with a boundary $\partial \Omega$ containing more than one point, there exists a conformal mapping $f\colon \mathbb{D} \to \Omega$, which is unique up to a Möbius transformation in \( \text{Aut}(\mathbb{D}) \).
\end{theorem}

The group of Möbius transformations that map the unit disk onto itself are called Möbius automorphisms of the unit disk, and denoted by \( \text{Aut}(\mathbb{D}) \). The general
form of a transformation in this group is given by:
\[ f(z) = e^{i\theta} \frac{z - \alpha}{1 - \overline{\alpha}z} \]
where:
\begin{itemize}
    \item \( \alpha \) is a complex number with \( |\alpha| < 1 \),
    \item \( \theta \) is a real number representing the angle of rotation,
    \item \( e^{i\theta} \) is a complex exponential indicating a rotation,
    \item \( \overline{\alpha} \) is the complex conjugate of \( \alpha \).
\end{itemize}

It should be noted that while a conformal mapping is only defined at the interior points of the domain, by the Carath\'eodory extension theorem (see e.g. \cite[Theorem 5.1.1]{Krantz2006}), a mapping between Jordan domains admits a homeomorphic extension to the boundary. This leads to a natural question known as the Gr\"otzsch problem. In fact, it is known that fixing images of three distinct boundary points is sufficient to make a conformal mapping between two Jordan domains unique. 

Let $\Omega$ be a Jordan domain in the complex plane $\mathbb{C}$, and let $z_1, z_2, z_3, z_4$ be four distinct positively oriented points on the boundary $\partial \Omega$. We call such domain together with the selected boundary points a (generalized) quadrilateral, and denote it by $Q = (\Omega; z_1, z_2, z_3, z_4)$.

The above considerations motivate the following definition \cite[p. 54, Definition 2.1.3]{PapamichaelStylianopoulos2010}:

\begin{definition}[Conformal modulus of a quadrilateral]
\label{def:conformal_modulus}
 If there exists a conformal mapping which transforms a quadrilateral $Q$ onto a rectangle $R_h = \big((0,1)\times (0,h); 1+i h, i h, 0, 1\big)$ so that the boundary points $z_1, z_2, z_3, z_4$ are mapped onto the vertices $1+i h, i h, 0, 1$, respectively, we call this number $h>0$ the (conformal) modulus of the quadrilateral $Q$, and denote it by
\begin{equation}
    \mathcal{M}(Q) = h.
\end{equation}
\end{definition}

In view of the above we note that the modulus $h$ is a unique characteristic that determines the conformal equivalence class of the quadrilateral.

We denote by $\widetilde{Q} = (\Omega; z_2, z_3, z_4, z_1)$ the so-called conjugate quadrilateral of $Q$. From the elementary properties of conformal moduli, one may observe the following basic identity (see e.g. \cite[p. 15]{Ahlfors2006}):

\begin{lemma}[Reciprocal identity]
Let $Q=(\Omega; z_1, z_2, z_3, z_4)$ and $\widetilde{Q} = (\Omega; z_2, z_3, z_4, z_1)$. Then
\begin{equation}
    \mathcal{M}(Q) \mathcal{M}(\widetilde{Q}) = 1.
    \label{eq:reciprocal_identity}
\end{equation}
\end{lemma}

This principle plays an important role in the study of geometric properties of quadrilaterals in the complex plane, as detailed in \cite[p.~15]{LehtoVirtanen1973} and \cite[pp.~53--54]{PapamichaelStylianopoulos2010}.

\subsection{Conformal moduli and Dirichlet--Neumann mixed boundary value problems}

Recall that the modulus of a quadrilateral $Q$ is closely related to a the Dirichlet--Neumann mixed boundary value problem(Dirichlet--Neumann problem, in short) of Laplace equation, see e.g. \cite[p.~431]{Henrici1986}.

\begin{definition}[Dirichlet--Neumann problem in general]
\label{def:general-DN-problem}
Consider a domain $\Omega$ where its boundary $\partial\Omega$ is a regular Jordan curve. The exterior unit normal $n$ is defined at almost every boundary point. Suppose that $\partial \Omega = \Gamma_{D} \cup \Gamma_{N}$, and  $\Gamma_{D}, \Gamma_{N}$ are composed of disjoint sets of regular Jordan arcs, intersecting at a finite number of points. Define $\beta_{D}$ and $\beta_{N}$ as continuous real-valued functions on $\Gamma_{D}$ and $\Gamma_{N}$, respectively. Let $\partial/\partial n$ denote differentiation to the direction of the exterior normal. The goal is to find a function $u$ such that:
\begin{enumerate}
    \item $u$ is continuous in $\overline{\Omega}$ and differentiable in $\Omega$.
    \item For every point $z \in \Gamma_{D}$, $u(z) = \beta_D(z)$.
    \item For each $z \in \Gamma_{N}$, the normal derivative $\frac{\partial}{\partial n} u(z) = \beta_N(z)$.
\end{enumerate}
\end{definition}




For quadrilateral $Q = (\Omega; z_1, z_2, z_3, z_4)$, we consider the boundary segments $\gamma_j$, for $j = 1, 2, 3, 4$, as the parts of $\partial \Omega$ connecting $z_1$ to $z_2$, $z_2$ to $z_3$, $z_3$ to $z_4$, and $z_4$ back to $z_1$, respectively. Let $\Gamma_{D} = \gamma_2\cup\gamma_4$,  $\Gamma_{N} = \gamma_1\cup\gamma_3$ and $u$ be the (unique) harmonic solution to the Dirichlet--Neumann mixed boundary value problem: 
\begin{equation}
\label{dir-neu}
\begin{cases}
    \Delta u(z) = 0 & \text{if } z \in \Omega, \\
    u(z) = 0 & \text{if } z \in \gamma_2, \\
    u(z) = 1 & \text{if } z \in \gamma_4, \\
    \frac{\partial}{\partial n} u(z) = 0 & \text{if } z \in \gamma_1 \cup \gamma_3.
\end{cases}
\end{equation}

Then by Ahlfors \cite[Theorem~4.5]{Ahlfors1973} and Papamichael and Stylianopoulos \cite[Theorem~2.3.3]{PapamichaelStylianopoulos2010}) 
\begin{equation}
\mathcal{M}(Q)=\iint_{\Omega}|\nabla u|^2 \, dx \, dy.
\end{equation}
The problem for the conjugate quadrilateral $\widetilde{Q}$ is called the conjugate Dirichlet--Neumann problem.

\section{Conjugate Function Method}\label{sec-conj}

Suppose that $Q$ is a quadrilateral, and $u$ is the harmonic solution to the Dirichlet--Neumann problem \eqref{dir-neu}. Let $v$ be the harmonic conjugate of $u$, normalized by the condition $v(\operatorname{Re} z_3, \operatorname{Im} z_3)=0$. Then the analytic function $f=u+iv$  maps $\Omega$ conformally onto the rectangle $R_h$, so that vertices are mapped onto $1+i h$, $i h$, $0$, and $1$. See \cite{HAKULA2013340}.

\begin{figure}[ht]
\centering
\includegraphics[max width=0.7\textwidth]{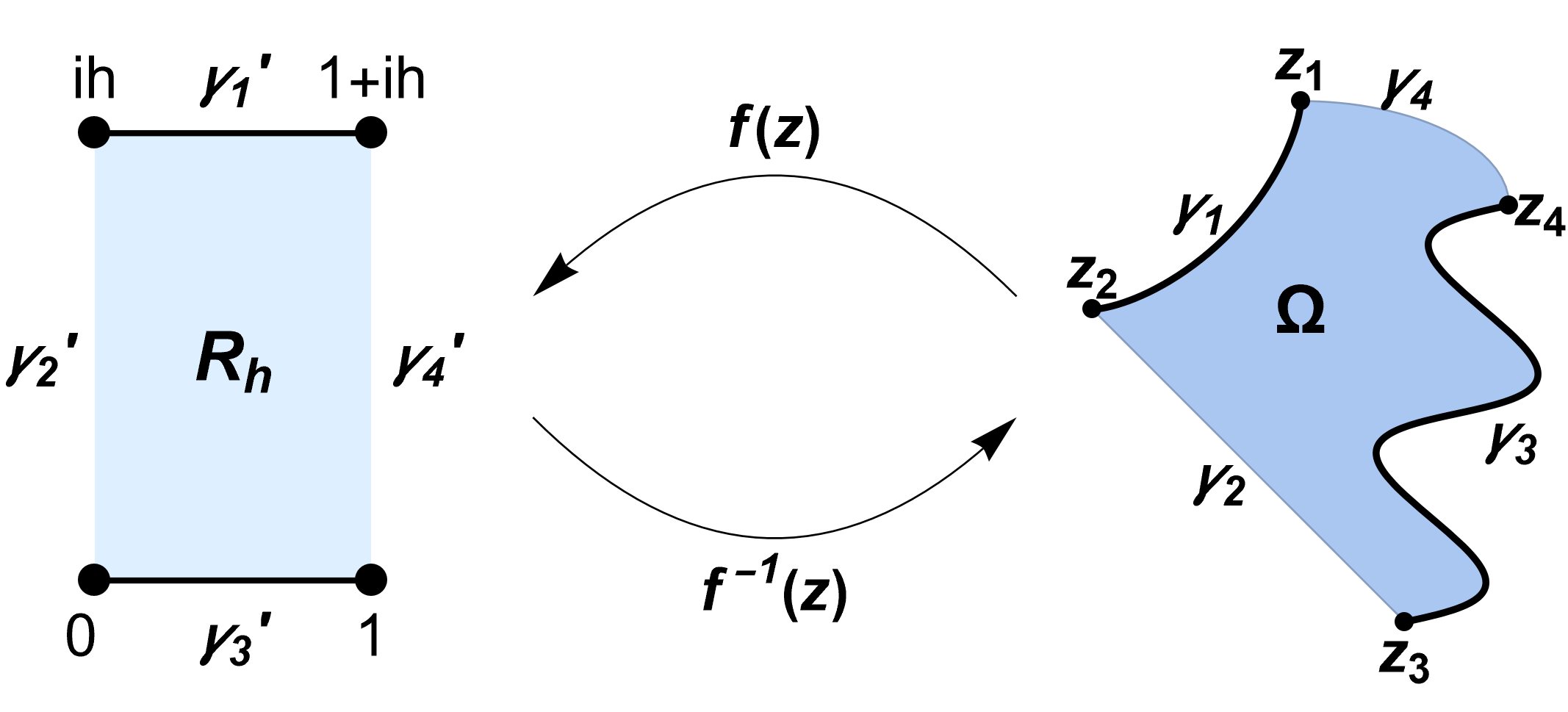}
\caption{Boundary value problem \eqref{dir-neu}. The Dirichlet boundary conditions are marked as thick lines. The Neumann boundary conditions
are thin lines.}
\label{fig:DirichletNeumannProblem}
\end{figure}

Now we may observe that  the canonical conformal mapping of a quadrilateral $Q = (\Omega; z_1, z_2, z_3, z_4)$ onto the rectangle $R_h$ with vertices at $1+i h$, $i h$, $0$, and $1$, can be obtained from solving the corresponding Dirichlet--Neumann problem $u$ and its conjugate problem. The following lemma (see \cite[Lemma 2.3]{HAKULA2013340}) shows how to construct the conjugate harmonic function $v$.

\begin{lemma}
Let $Q$ be a quadrilateral with modulus $h$, and suppose $u$ solves the Dirichlet--Neumann problem \eqref{dir-neu}. If $v$ is a harmonic function conjugate to $u$, satisfying $v(\operatorname{Re} z_3, \operatorname{Im} z_3)=0$, and $\tilde{u}$ represents the harmonic solution for the conjugate quadrilateral $\widetilde{Q}$, then $v = h \tilde{u}$.
\end{lemma}

In particular, if $u$ and $\tilde{u}$ are known, we may derive the value of the modulus $h$ through the Cauchy--Riemann equations. Once $h$ is known, we get the conformal mapping. The conjugate function method for constructing a conformal mapping is as follows:

\subsection*{Steps for Finding Conformal Mapping}

\begin{enumerate}
    \item Find the solution \( u \) to the Dirichlet--Neumann problem \eqref{dir-neu}.
    \item Find the solution \( \tilde u \) to the Dirichlet--Neumann problem associated with \( \widetilde{Q} \).
    \item Pick a suitable point inside the domain (usually near the center of the domain) and evaluate suitable pairs of partial derivatives \( \{ u_x, \tilde u_y \} \) and \( \{ \tilde u_x,  u_y \} \), where $z=x+iy$. If \( f = u + i h \tilde u \) is a conformal mapping, it satisfies the Cauchy--Riemann equations. Therefore we can compute $h$ as
\[
 h = \frac{u_x}{\tilde u_y} \text{      or      }  h = - \frac{u_y}{\tilde u_x}
\]
The vertices \( \{z_1, z_2, z_3, z_4\} \) are mapped onto the corners \( \{1+i h, i h, 0, 1\} \).
\end{enumerate}

The partial derivatives of $u$ and $\tilde u$ can be estimated using the five-point formula \cite[p. 883]{AbramowitzStegun}
:
\[
u_x \approx \frac{-u(x+2\Delta x, y) + 8u(x+\Delta x, y) - 8u(x-\Delta x, y) + u(x-2\Delta x, y)}{12\Delta x}.
\]
The error term for this approximation depends on the fourth derivative of $u(x, y)$ and is given by:
\[
\text{Error} = -\frac{(\Delta x)^4}{30} \frac{\partial^4 u}{\partial x^4}(\xi, y),
\]
where $\xi$ is a point in the interval $[x-2\Delta x, x+2\Delta x]$.

It is clear that solving the Dirichlet--Neumann problem \eqref{dir-neu} is the key step. Section \ref{sec-walk} will introduce a stochastic approach to solve it numerically through simulation, which we call the Reflected Walk-on-Spheres (RWoS) method.

\section{Theory behind the WoS method}\label{sec-walk}

The Walk-on-Spheres (WoS) method is a Monte Carlo algorithm used for solving the Dirichlet problem for the Laplace equation. 

In this section, we first introduce the Brownian motion and state its invariance with respect to conformal mappings. Then we show the relation between Dirichlet--Neumann boundary value problems and Brownian motions, which is the foundation for our Reflected Walk-on-Spheres method.


\subsection{Brownian motion}

\begin{definition}
A \emph{one-dimensional Brownian motion} is a stochastic process \( \{ B(t), t \geq 0 \} \) characterized by:
\begin{enumerate}
    \item The starting condition \( B(0) = x \) for a given real number \( x \).
    \item Independence of increments: For any non-overlapping intervals \( [s_1, t_1] \) and \( [s_2, t_2] \), the increments \( B(t_1) - B(s_1) \) and \( B(t_2) - B(s_2) \) are independent.
    \item For each time interval \( [t, t+h] \), the increment \( B(t+h) - B(t) \) is normally distributed with a mean of 0 and a variance of \( h \).
    \item Continuity in probability: The path \( t \mapsto B(t) \) is continuous with probability 1.
\end{enumerate}
The case where \( B(0) = 0 \) defines the standard Brownian motion.

A \emph{multi-dimensional Brownian motion} in \( d \) dimensions, denoted as \( \{ \mathbf{B}(t) = (B_1(t), \ldots, B_d(t)), t \geq 0 \} \), consists of \( d \) independent one-dimensional Brownian motions.  In particular, the complex Brownian motion is $B= B_1 + iB_2$, where $B_1$ and $B_2$ are two independent real Brownian motions.
 \end{definition}


Recall the following  definition of the harmonic measure. For details we refer to \cite[Chapter I]{GarnettMarshall2005}.

\begin{definition}[Harmonic measure]
For the upper half-plane \( \mathbb{H} \) and a point \( z \in \mathbb{H} \), the harmonic measure \( \omega(z, E, \mathbb{H}) \) of a set \( E \)  on the real line at the point $z$ is defined as:
\begin{equation}
\omega(z, E, \mathbb{H}) = \sum_{j=1}^{n} \frac{\theta_j}{\pi}
\end{equation}
where \( E \) is a finite union of open intervals on the real line, and each \( \theta_j \) is given by:
\begin{equation*}
 \theta_j = \text{arg}\left(\frac{z - b_j}{z - a_j}\right) 
\end{equation*}
Here, \( (a_j, b_j) \) are the intervals so that \( E \) is their union, and \( \text{arg} \) is the argument of a complex number.

The harmonic measure has the following properties:
\begin{itemize}
\item[(i)] It is always between 0 and 1 for \( z \in \mathbb{H} \).
\item[(ii)] It approaches 1 as \( z \) approaches \( E \).
\item[(iii)] It approaches 0 as \( z \) approaches the complement of \( E \) in the extended real line.
\end{itemize}
This measure is uniquely determined and is actually the unique harmonic function that satisfies the above properties. Its uniqueness is a consequence of Lindelöf's maximum principle.

For the unit disk \( \mathbb{D} \) , let \( E \) be a finite union of open arcs on \( \partial\mathbb{D} \). The harmonic measure of \( E \) at \( z \) in \( \mathbb{D} \) is defined to be
\begin{equation}
\label{eq: half plane harmonic measure}
\omega(z, E, \mathbb{D}) = \omega(f(z), f(E), \mathbb{H})    
\end{equation}
where \( f \) is any conformal mapping from \( \mathbb{D} \) onto \( \mathbb{H} \) such that $f(\partial\Omega)$ is the real line. This harmonic function also satisfies the three conditions similar to (i), (ii), and (iii). In addition by Lindelöf's maximum principle this definition is independent of how \( f \) is chosen. Moreover, by changing the variable $\phi(z) = i (1 + z) / (1 - z)$ we get
\begin{equation}
\omega(z, E, \mathbb{D}) = \int_{E} \frac{1 - |z|^2}{|e^{i\theta} - z|^2} \frac{d\theta}{2\pi}
\end{equation}

Let \( \Omega \) be a simply connected domain and \( \partial\Omega \) is a Jordan curve.  Let $\phi$ be a conformal mapping from the unit disk $\mathbb{D}$ onto $\Omega$, $z$ is a point in $\Omega$ and $w=\phi^{-1}(z)$. For any Borel set $E\subset \partial\Omega$ the harmonic measure of $E$ with respect to $z$ in $\Omega$ is defined as:
\begin{equation}
\label{eq: harmonic measure for omega}
\omega(z, E, \Omega) = \omega(\phi^{-1}(z), \phi^{-1}(E), \mathbb{D}) = \int_{\phi^{-1}(E)} \frac{1 - |w|^2}{|e^{i\theta} - w|^2} \frac{d\theta}{2\pi}
\end{equation}

\end{definition}

\begin{remark}
The harmonic measure in \eqref{eq: harmonic measure for omega} does not depend on the choice of $\phi$. If there is another function $\psi$ that maps $\mathbb{D}$ onto $\Omega$, then by Riemann mapping theorem, $\phi^{-1}$ and $\psi^{-1}$ differ only by a Möbius automorphism $T$ of $\mathbb{D}$, i.e. $\psi^{-1} = T \circ 
 \phi^{-1} $. Therefore by \eqref{eq: half plane harmonic measure},
\begin{equation}
\begin{split}
\omega(\psi^{-1}(z), \psi^{-1}(E), \mathbb{D}) &= \omega\big( T\circ (\phi^{-1}(z)), T\circ (\phi^{-1}(E)), \mathbb{D}\big) \\
&= \omega\big( (f \circ T)\circ (\phi^{-1}(z)), (f \circ T)\circ (\phi^{-1}(E)), \mathbb{H}\big) \\
&= \omega\big(\phi^{-1}(z), \phi^{-1}(E), \mathbb{D}\big),
\end{split}
\end{equation}
because $f \circ T$ is a conformal mapping.

Furthermore, if we choose the conformal mapping $\phi$ from $\Omega$ to $\mathbb{D}$ such that $\phi(z_0) = 0$ for $z_0 \in \Omega$, the harmonic measure of $E \subset \partial\Omega$ at $z_0$ in $\Omega$ becomes 
\begin{equation}
\omega(0, \phi^{-1}(E), \mathbb{D}) = \int_{\phi^{-1}(E)} \frac{1 - |0|^2}{|e^{i\theta} - 0|^2} \frac{d\theta}{2\pi}  = |\phi^{-1}(E)|
\end{equation}\
where $|\cdot|$ denotes the Euclidean length of an arc on the unit circle.

In this case the interpretation of harmonic measure becomes closely tied to probability theory. Specifically, it represents the probability that a Brownian motion starting from the origin, exits the disk through the set $\phi^{-1}(E)$. Due to the symmetry of the disk, this probability is uniform and is equal to $|\phi^{-1}(E)|/(2\pi)$. This illustrates the probabilistic interpretation of the harmonic measure.
\end{remark}


\subsection{Conformal invariance of the Brownian motion and the harmonic measure}

Brownian motion in the plane, or the complex Brownian motion, exhibits a conformal invariance, see e.g. \cite[Theorem 2.2]{Lawler2005}. Indeed, let $f$ be a conformal mapping from domain $\Omega$ onto $\Omega'$. If $B$ is a complex Brownian motion on $\Omega$ then $f(B)$ admits the representation
$$
f(B_t) = \widetilde B_{\int_0^t |f'(B_u)|^2 d u},
$$
where $\widetilde B$ is another complex Brownian motion.
Thus the image trajectory of a Brownian motion under conformal mapping is still a trajectory of Brownian motion.

Thus if we choose the conformal mapping $\phi$ from $\Omega$ to $\mathbb{D}$ such that $\phi(z_0) = 0$, when a Brownian particle starting at that point $z_0$ in $\Omega$ exits the domain through $E$, we would observe in $\mathbb{D}$ that the Brownian particle starting at 0 and exits the domain through $\phi^{-1}(E)$. So we can conclude the probability of $z$ hitting $E$ is precisely the harmonic measure $\omega(z, E, \Omega) = |\phi^{-1}(E)| $. More generally speaking, the probability of the Brownian motion hitting a particular part of the boundary in the conformally transformed domain is the same as in the original domain \cite{Lawler2005}. 


The next definition introduces a powerful function for dealing with reflection at Neumann boundaries in RWoS algorithm.

\begin{definition}[Anti-conformal Mapping]
An \emph{anti-conformal mapping} is a continuous function \( w = f(z) \) from a neighborhood of \( z_0 \) in the complex \( z \)-plane to a neighborhood of \( w_0 \) in the complex \( w \)-plane. It preserves angles but reverses orientation.
\end{definition}

 The reversal of orientation in anti-conformal mappings changes the trajectories of the Brownian paths by mirroring them. Consequently, a Brownian trajectory remains a Brownian trajectory under anti-conformal maps.

\begin{proposition}
The trajectory of a Brownian motion under an anti-conformal map is again a trajectory of a Brownian motion in the image domain. The same holds for the harmonic measure. Indeed, \( \omega(X, \Omega, z) = \omega(f(X), f(\Omega), f(z)) \) for any subset \( X \subset \partial \Omega \) and anti-conformal mapping \( f \).
\end{proposition}

\subsection{Connection between Dirichlet--Neumann boundary value problems and Brownian motion}

The connection between Dirichlet boundary value problems and Brownian motion is due to Kakutani in 1944 \cite{Kakutani1944}.
Kakutani showed that the solution to the Dirichlet problem can be represented in terms of Brownian motion. Specifically, if a Brownian motion starts from a point inside \( \Omega \), the value of the harmonic function \( u \) at that point is the expected value of the dirichlet boundary function \( \beta \) evaluated at the point where the Brownian motion first exits the domain \( \Omega \).
Indeed, if \( B \) is a Brownian motion starting at point \( z = x + iy \) inside \( \Omega \), and \( \tau \) is the first exit time from \( \Omega \), then the solution \( u(z) \) to the Dirichlet problem is given by:
\[ u(z) = \mathbb{E}_z[\beta(B(\tau))] \]
Here, \( \mathbb{E}_z[f(B(\tau))] \) is the expected value of \( f \) at the exit point of the Brownian motion from \( \Omega \) starting from $z = x + iy\in\Omega$.

In the plane, the Kakutani connection has been extended for the Dirichlet--Neumann problem, see e.g. Sylvain Maire and Etienne Tanré \cite{MaireTanre2013}. When the normal derivative at the Neumann Boundary is zero, the stochastic representation is simplified enough to simulate through WoS method.  In this case, Brownian motion is reflected about the normal line at the point where it hits the Neumann boundary (a process known as specular reflection). Following this reflection, the (reflected) Brownian motion is halted at the first time it hits the Dirichlet boundary. For the proof, we refer to \cite{MaireTanre2013}.

\begin{proposition}[Stochastic Representation]
 Let \( \Omega \) be a bounded domain in \( \mathbb{C} \) with a regular boundary \( \Gamma = \Gamma_D \cup \Gamma_N \)  where $\Gamma_{D}$ and $\Gamma_{N}$ are both unions of regular Jordan arcs such that $\Gamma_{D}\cap\Gamma_{N}$ is finite. If \( \beta_{N}=0 \) and \( \beta_{D} \) is a continuous function on \( \Gamma_{D} \).  Then, For a point \( z = x + iy \) in \( \Omega \), the harmonic solution \( u \) to the Dirichlet--Neumann problem is given by \( u(z) = \mathbb{E}_z\left[\beta_{D}\left(B_\tau\right)\right] \), where $B$ is a $2$-dimensional Brownian motion, the expected value is taken conditionally on $\left\{B_0=z = x+iy\right\}$, and $\tau$ is the first exit time from $\Gamma_{D}$, which is also the first-hitting time of $\Gamma_{D}$. When the particle touches the Neumann boundary $\Gamma_{N}$, it is reflected symmetrically back into the domain. 
\end{proposition}

The Dirichlet--Neumann problem in \eqref{dir-neu} is the case when $\Gamma_{D}=\gamma_{2}\cup\gamma_{4}$, $\Gamma_{N}=\gamma_{1}\cup\gamma_{3}$, $\beta_{D}=1$ on $\gamma_{4}$ and $\beta_{D}=0$ on $\gamma_{2}$.

To approximate a solution to the Dirichlet--Neumann problem using Brownian motion, one can simulate multiple independent Brownian paths starting from a point inside the domain. By the Law of Large Numbers, the average value of a function evaluated at the exit points of these paths will converge to the expected value, providing an approximation of the solution. Mathematically, this is expressed as:
\[
\mathbb{E}_z\left[\beta_{D}\left(B_\tau\right)\right] \sim \frac{1}{k} \sum_{i=1}^k \beta_{D}\left(B_{\tau_i}^i\right)
\]
where \( k \) is the number of simulated paths, and \( B_{\tau_i}^i \) is the exit point of the \( i \)-th path.

\section{Reflected Walk-On-Spheres Algorithm}\label{sec-refl}

An obvious method for simulating Brownian motion involve discretizing its path into small time intervals and computing the movement in each interval. This can be computationally intensive, especially for fine time scale discretizations. The WoS method, on the other hand, uses a different approach by jumping directly from one point to another within a sphere, significantly reducing the number of steps needed to simulate the Brownian path.
The Walk-on-Spheres (WoS) method \cite{Muller} exploits the rotational invariance of Brownian motion: the exit point of a Brownian trajectory from a sphere is uniformly distributed over the sphere's surface. The WoS method starts at an initial point in the domain and repeatedly draws the largest possible sphere within the domain, selecting the next point uniformly from the sphere's surface until the boundary is hit. This process is repeated to build a sequence of points that simulate the Brownian motion on the spheres. Versions of WoS method for eigenvalue type equations have been investigated in \cite{YRS1,YRS2}.

\subsection{Sampling the Exit Point}
The WoS algorithm samples the first exit point of a Brownian motion from a domain.  Starting with \( z^{(0)} = z \), it iteratively draws the largest possible sphere \( \mathcal{S}_0 \) within the domain and determines the exit point \( z^{(1)} \) until the Brownian motion approaches the boundary.  One may think the process converges to the boundary very fast. Yet, practically, the algorithm would require an infinite number of steps to stop \cite{MascagniHwang2003} since the imaginary boundary has no thickness and the closer the Brownian particle is to the boundary, the smaller the radius of the sphere is. To address this in computational implementations, the process is halted as soon as it approaches close enough to the boundary, typically within an \( \varepsilon \)-shell, and the point of exit is projected onto the boundary. This approach is often referred to as the introduction of an 
$\epsilon$-shell or  $\epsilon$-layer.

\begin{figure}[h]
    \centering
    \includegraphics[max width=\textwidth]{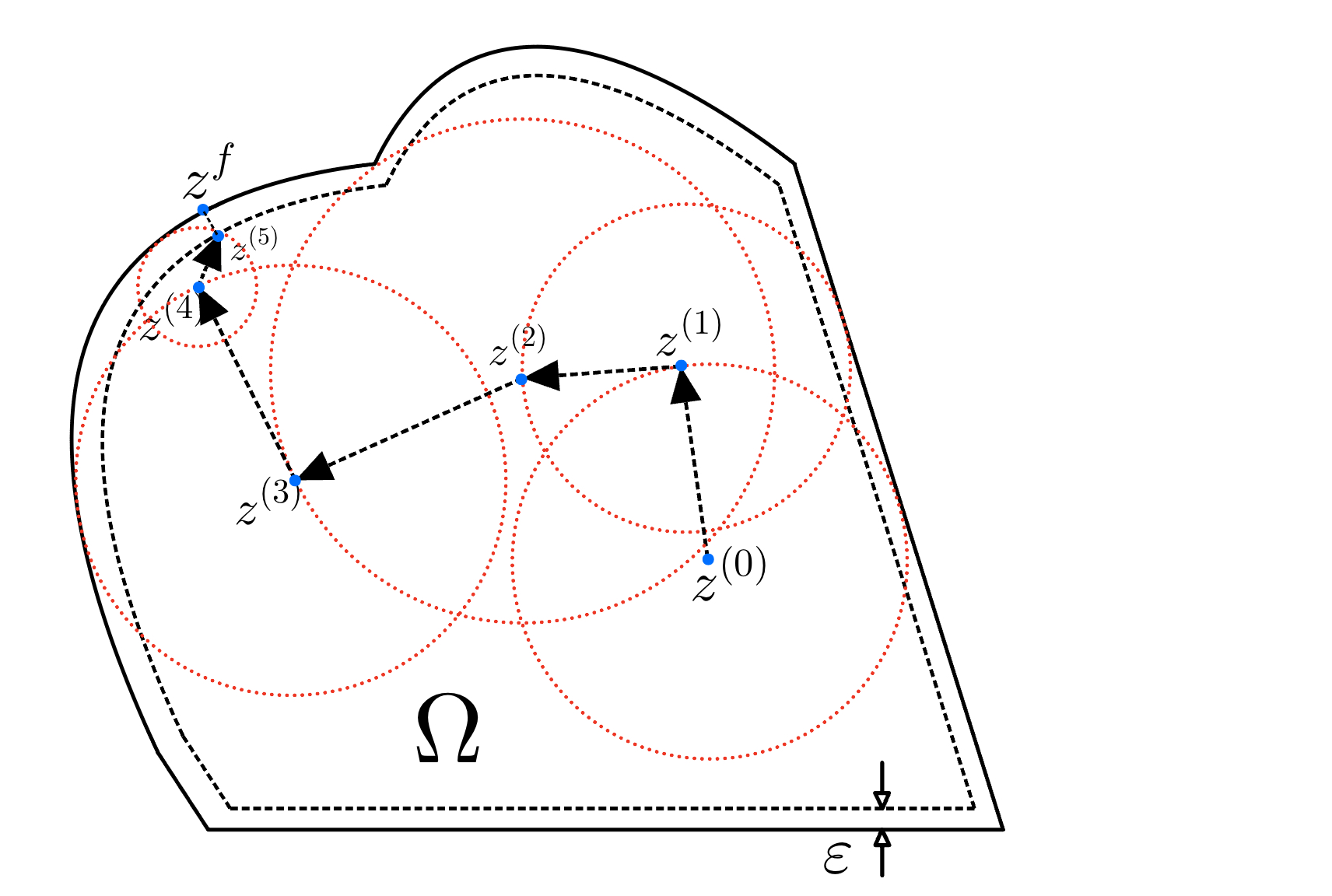}
    \caption{An illustration of the Walk-on-spheres algorithm with an $\varepsilon$-shell.}
    \label{fig:wos}
\end{figure}

\subsection{Reflection at Neumann Boundaries}
In the Dirichlet--Neumann problem the Brownian trajectory is absorbed in the Dirichlet boundary and reflected on the zero-valued Neumann boundary. The challenge in applying the WoS algorithm in solving the 2d Dirichlet--Neumann problem is to accurately reflect the Brownian motion at Neumann boundaries. One idea is that, when the particle is near the Neumann boundary, stop walking on spheres and use the primitive idea of simulating fine time mesh to compute the reflection angle about the normal lines. 
But this is not only very costly in computation, but also it brings in an error dependent on step size in the simulation. Therefore, the idea of fine time mesh to simulate reflection is not very desirable. 
Fortunately, for Neumann boundaries that are straight segments, we can think of them as mirrors and walk on spheres through the mirrors as if there is no barrier between our real world and the reflected imaginary world. As soon as the Brownian particle stepped into the reflected domain, we bring it back by reflecting it about the straight Neumann boundary, just as if it had bounced off the mirror (maybe more than once). In this way, we have walked on a sphere not bounded by the Neumann boundary. Although we don't know the exact reflection process the particle went through, we have simulated the destination point at a macroscopic level. And then we set this destination point as the center for the next sphere walk.


\subsection*{Constructing the anti-conformal map for reflection over curved boundary}

In addition of straight lines, we also consider curved Neumann boudaries.  The reflection with respect to curved boundaries are handled by using anti-conformal mappings.

In the following $N\subseteq\Gamma_{N}\subseteq\partial\Omega$ is a curved Neumann boundary.

\begin{enumerate}
    \item Find the conformal mapping $f$ that transform $N$ into a straight line segment $I$.
    \item Find the reflection function $R$ about the line segment $I$.
    \item Then the composition function $g = f^{-1} \circ R \circ f$ is the desired anti-conformal mapping which mirrors $\Omega$ into the reflected domain $g(\Omega)$. Since only the set $N$ remains fixed in the reflection, it would be expected $\Omega \cap g(\Omega) = N$.  
\end{enumerate}

In addition, $g = g^{-1}$ since
\begin{eqnarray*}
g \circ g &=& (f^{-1} \circ R \circ f) \circ (f^{-1} \circ R \circ f)  \\ 
&=& f^{-1} \circ R \circ (f \circ f^{-1}) \circ R \circ f  \\
&=& f^{-1}\circ ( R \circ R) \circ f\\
&=& f^{-1} \circ f = \text{Identity}.
\end{eqnarray*}
Therefore, $g(g(\Omega)) = \Omega$.

The remaining question is how to find the conformal mapping $f$ from curved boundary to line segment. This is a difficult question and there is no general answer yet. However, this issue is mitigated for circular-arc boundaries for which we have the well-known Möbius transformations that can map circles into lines.

If we use the composition mentioned above, for circular-arc boundaries with radius \( R \) centered at \( (x_1, y_1) \), the mapping 
\begin{equation}
\label{eq:anti-conformal_mapping_for_reflection}
(x,y)\mapsto \left(\frac{R^2 \left(x-x_1\right)}{\left(x-x_1\right){}^2+\left(y-y_1\right){}^2}+x_1,\frac{R^2 \left(y-y_1\right)}{\left(x-x_1\right){}^2+\left(y-y_1\right){}^2}+y_1\right)
\end{equation}
is the anti-conformal reflection. 
This mapping is an inversion in $S((x_1,y_1),R)$ and hence preserves the class of generalized circles\cite[Theorem~3.2.1]{Beardon1983}.
Thus the reflections at circular-arc Neumann boundary are manageable.

\begin{remark}
Note that the particle cannot be reflected somewhere outside the original domain $\Omega$. Also note that some line segments or circular arcs may adjoin each other to form a consecutive Neumann boundary. In this case it is necessary to break it into smaller tractable pieces.  Therefore when calculating each radius of next walking sphere, we should set it to be the minimal distance between the current position $z^{(n)}$ and $\Gamma_{D} \cup S \cup \left( \bigcup_{i} g_i\left(\partial\Omega \setminus N_i\right) \right)$, where $S$ is the splitting points of the consecutively adjoining Neumann boundaries, $g_i$ is the anti-conformal mapping function for reflecting about the $i$th Neumann boundary $N_i$ among all the Neumann boundaries $\Gamma_N \subseteq \partial\Omega$. In this way, it ensures the Brownian trajectories stay in the original domain. In practice, the distance calculation $d \left( z^{(n)}, \left( \Gamma_{D} \cup S \cup \left( \bigcup_{i} g_i(\partial\Omega \setminus N_i) \right) \right) \right)$ can often be optimized by not considering some boundaries due to the intrinsic geometry of the original domain $\Omega$.
    
\end{remark} 

\begin{algorithm}[H]
\caption{Reflected Walk-on-spheres (RWoS) method for sampling reflected Brownian trajectory exit points in polygonal or circular-arc domains}
\label{alg:walk_on_spheres}
\begin{algorithmic}[1]
\Require $\varepsilon > 0$, $\Omega$ is the domain in a plane. $\partial\Omega = \Gamma_D \cup \Gamma_N$ where $\Gamma_{D}$ is the set of Dirichlet boundaries and $\Gamma_{N}$ is the set of zero-normal-derivative Neumann boundaries. 
\State Split the $\Gamma_N$ into circular arcs or line segments $\bigcup_{i} N_i$, let $S$ be their splitting points. 
Let $g_i$ be the anti-conformal mapping function for reflecting about the $i$th Neumann boundary $N_i$ among all the Neumann boundaries $\bigcup_{i} N_i = \Gamma_N$.
\State Initialize: $z^{(0)} = z$
\While{$d(z^{(n)}, \Gamma_D) > \varepsilon$}
    \State  Set $r_n = d \left( z^{(n)}, \left( \Gamma_{D} \cup S \cup \left( \bigcup_{i} g_i(\partial\Omega \setminus N_i) \right) \right) \right)$. 
    \State Sample $\gamma_n$ uniformly distributed over the unit sphere.
    \State Set $z^{(n+1)} := z^{(n)} + r_n \gamma_n$.
    \If{$z^{(n+1)} \notin \Omega$}
        \If{Neumann boundary is a line segment}
            \State Reflect $z^{(n+1)}$ about the line back to $\Omega$.
        \Else \Comment{(Neumann boundary is a circular arc)}
            \State Reflect $z^{(n+1)}$ by function~\ref{eq:anti-conformal_mapping_for_reflection}.
        \EndIf
    \EndIf
\EndWhile
\State Set $z_f := p_{\Gamma}(z^{(n)})$, the orthogonal projection on $\Gamma$.
\State \Return $z_f$
\end{algorithmic}
\end{algorithm}

\begin{figure}[h]
    \centering
    \includegraphics[max width=\textwidth]{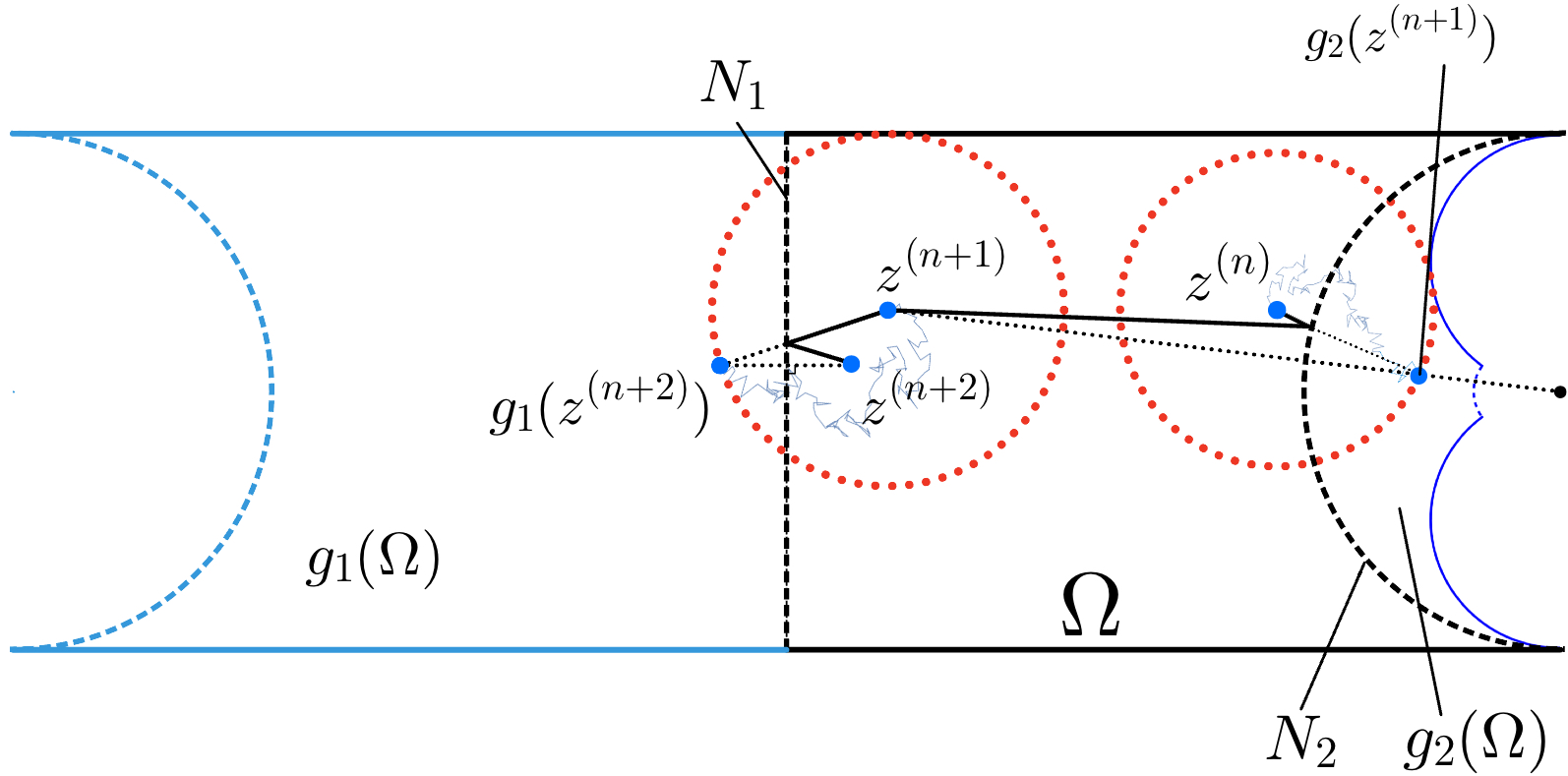}
    \caption{An illustration of the Reflected Walk-on-spheres algorithm. The original domain $\Omega$ is the region enclosed by a rectangle and a half circle. The left side of the rectangle $N_1$ and the half circle are the Neumann boundaries. $g_1, g_2$ are the anti-conformal mappings that reflects $\Omega$.}
    \label{fig:wos}
\end{figure}



\clearpage

\section{Examples}\label{sec-ex}

In this section, we give simulation examples and results to demonstrate our method is an effective one in numerically solving mixed boundary Laplace equations and constructing conformal mappings.

For the canonical mapping $f = u + vi$ of quadrilaterals, we plot the equipotential curves of $u$ and $v$. We also give the computed values of their moduli and compare them to the reference values. These examples have been computed using a straightforward implementation of the algorithm on C programming language. To compare these results with the deterministic $hp$-finite element method, the readers should refer to \cite{HakulaRasilaVuorinen2011}.

In addition to the 2d examples, we our RWoS method can be extended to solve Laplace equation in polyhedral domains with zero Neumann boundary conditions in the three-dimensional space. This kind of problems arises in engineering, e.g., heat transfer problems. 

\subsection{Rectangle}

We first test our method on the simplest canonical domain, which is the rectangle $\{z \in \mathbb{C}: 0<\operatorname{Re}(z)<h, 0<\operatorname{Im}(z)<1\}$. By Definition \ref{def:conformal_modulus}, the modulus of the normalized rectangle is its height $h$.

\begin{figure}[H]
\centering
\includegraphics[width=6cm]{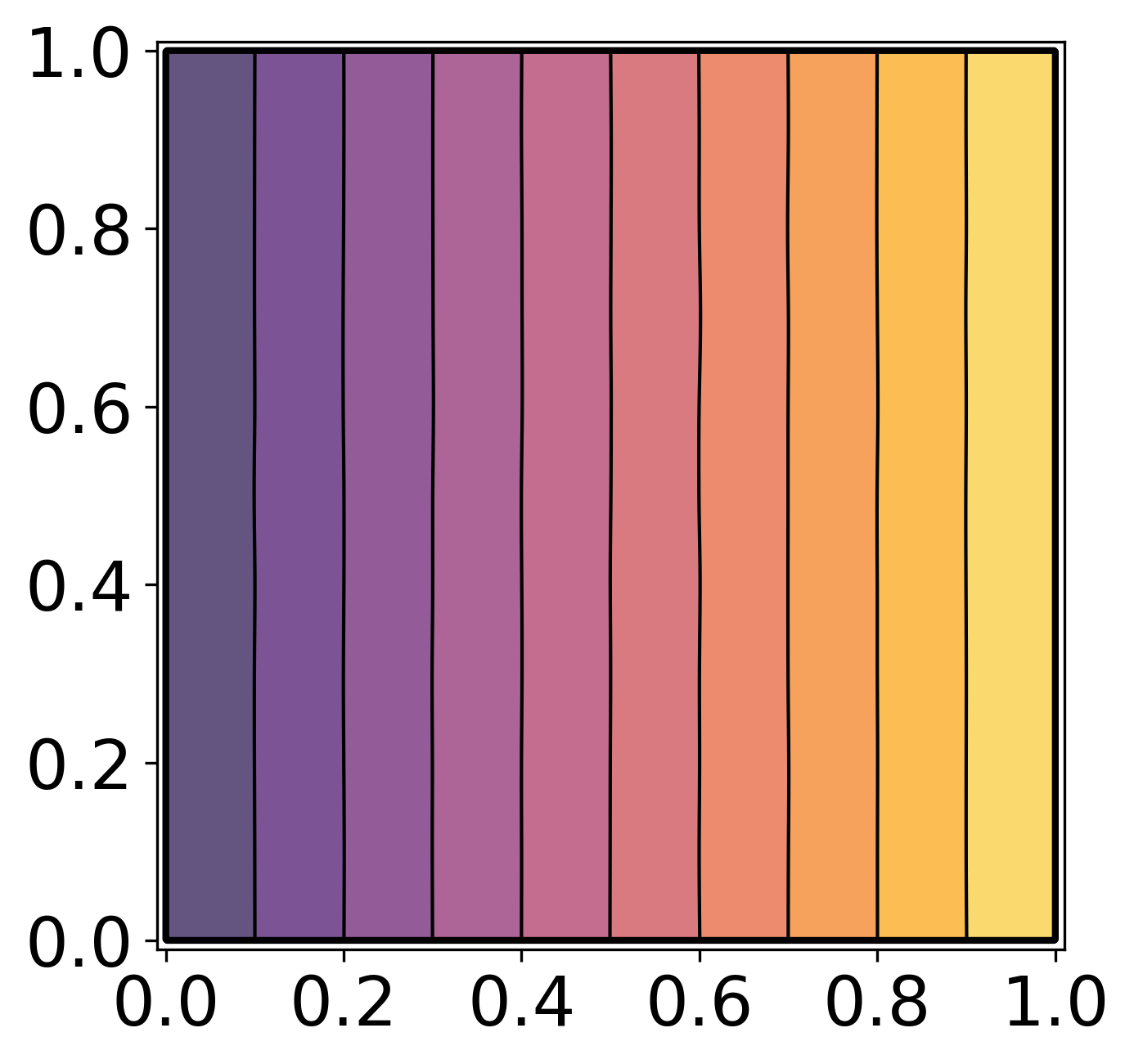}
\includegraphics[width=6cm]{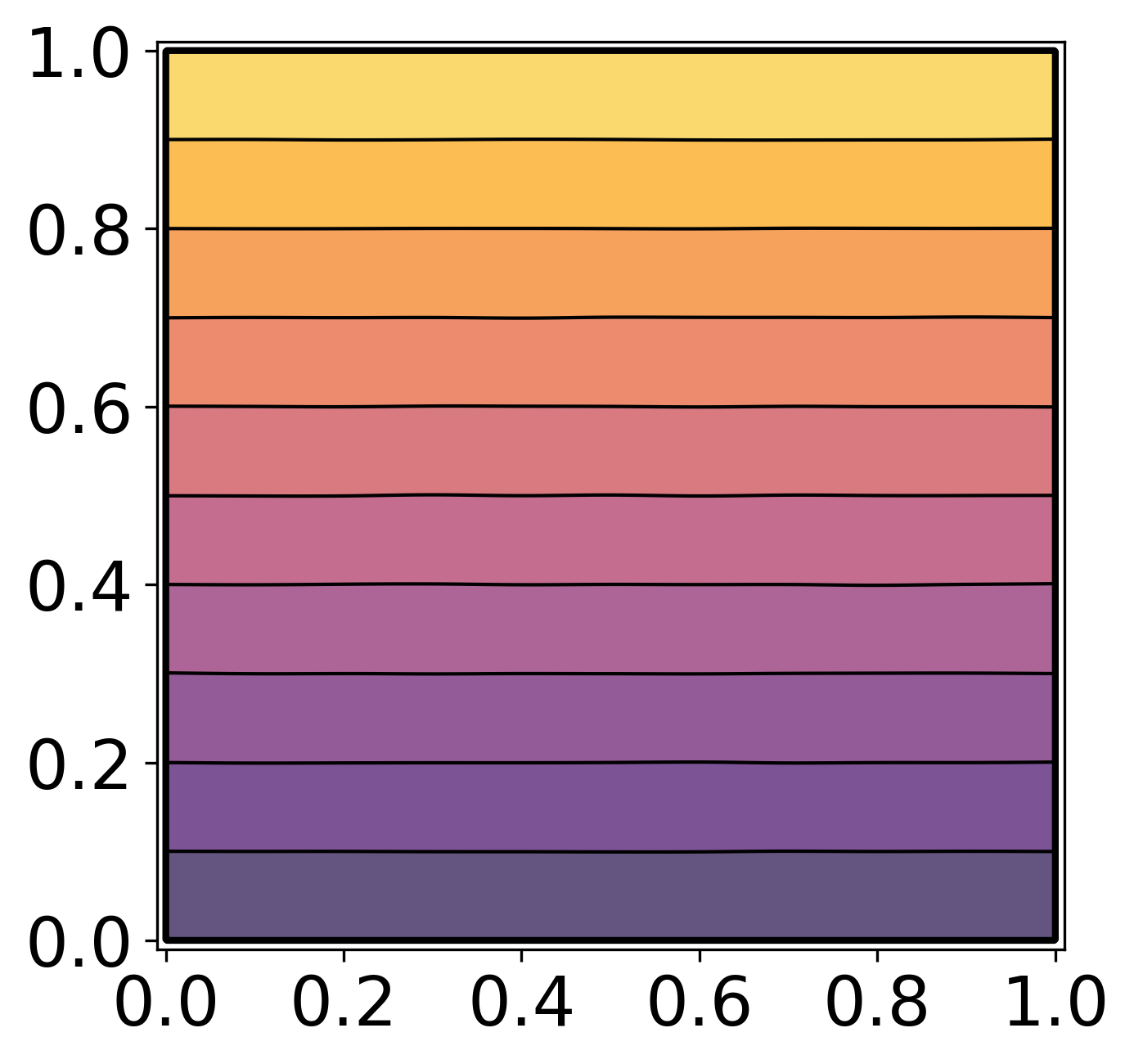}
\caption{The rectangle with vertices \((0, 0)\), \((1, 0)\), \((1, h)\), \((0, h)\).}
\label{fig:Rectangle_potential}
\end{figure}

\begin{table}[!htbp]
\centering

\begin{tabular}{|c|c|c|c|}
\hline
Height(Modulus) h & Computed value & Error \\
\hline
0.600000 & 0.600292 & $2.92 \cdot 10^{-4} $ \\
0.800000 & 0.799714 & $2.86 \cdot 10^{-4}$ \\
1.000000 & 1.000458 & $4.58 \cdot 10^{-4}$ \\
1.200000 & 1.200488 & $4.88 \cdot 10^{-4}$ \\
1.400000 & 1.400667 & $6.67 \cdot 10^{-4}$ \\
\hline
\end{tabular}
\caption{Moduli of rectangle with vertices \((0, 0)\), \((1, 0)\), \((1, h)\), \((0, h)\).}
\label{table:Rectangle_moduli}
\end{table}

\begin{figure}[H]
\centering
\includegraphics[width=0.5\textwidth]{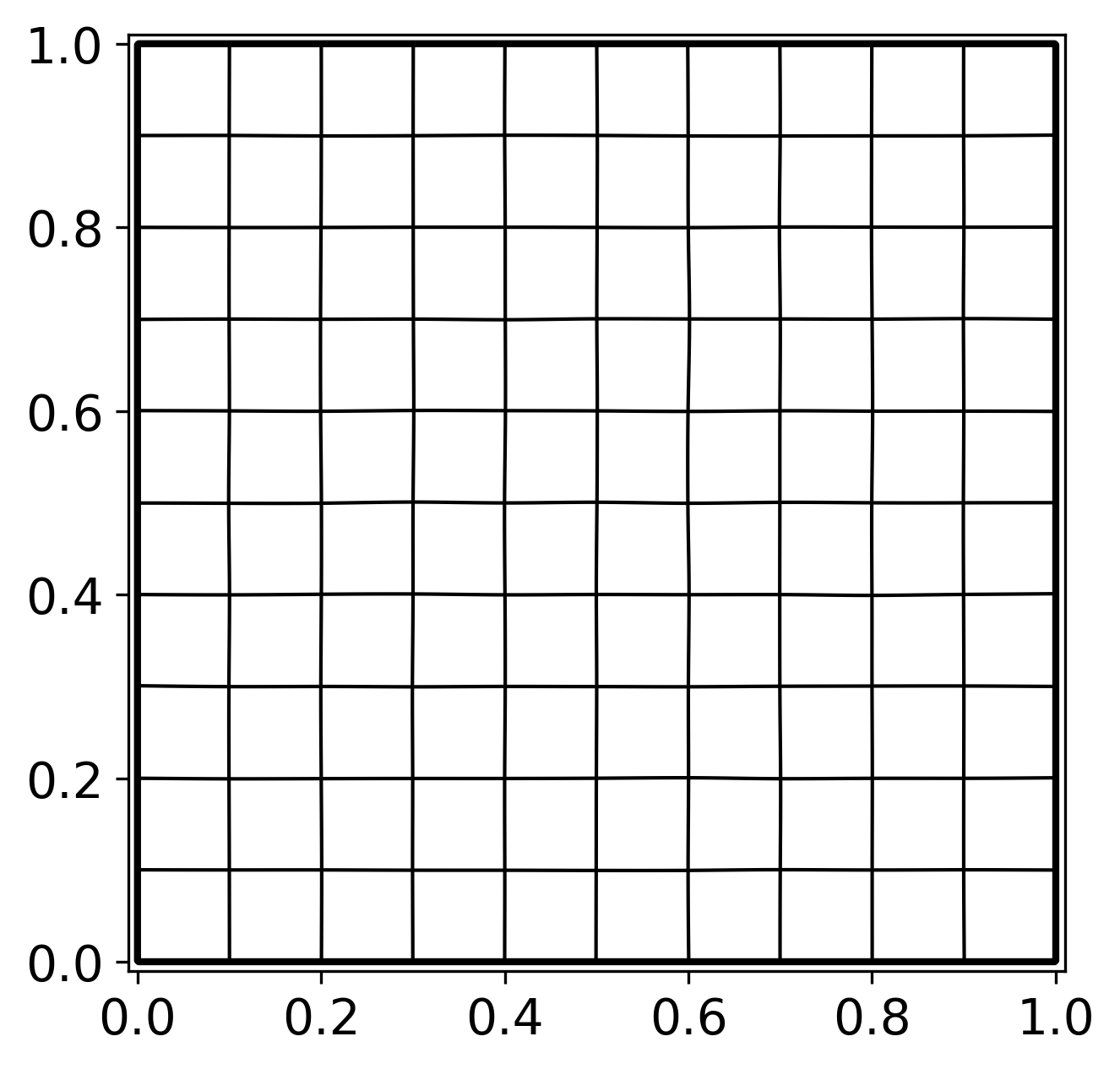}
\caption{The canonical conformal mesh for rectangle.}
\label{fig:Rectangle_mesh}
\end{figure}

\subsection{L-shaped region}

The L-shaped region
\[
\begin{gathered}
L(a, b, c, d)=L_1 \cup L_2, \quad L_1=\{z \in \mathbb{C}: 0<\operatorname{Re}(z)<a, 0<\operatorname{Im}(z)<b\}, \\
L_2=\{z \in \mathbb{C}: 0<\operatorname{Re}(z)<d, 0<\operatorname{Im}(z)<c\}, 0<d<a, 0<b<c
\end{gathered}
\]
is commonly used in various computational analyses. It has been explored by Gaier \cite{Gaier1995}, who provided computed moduli values for such domains. We compare our results with his.

\begin{figure}[H]
\centering
\includegraphics[width=6.2cm]{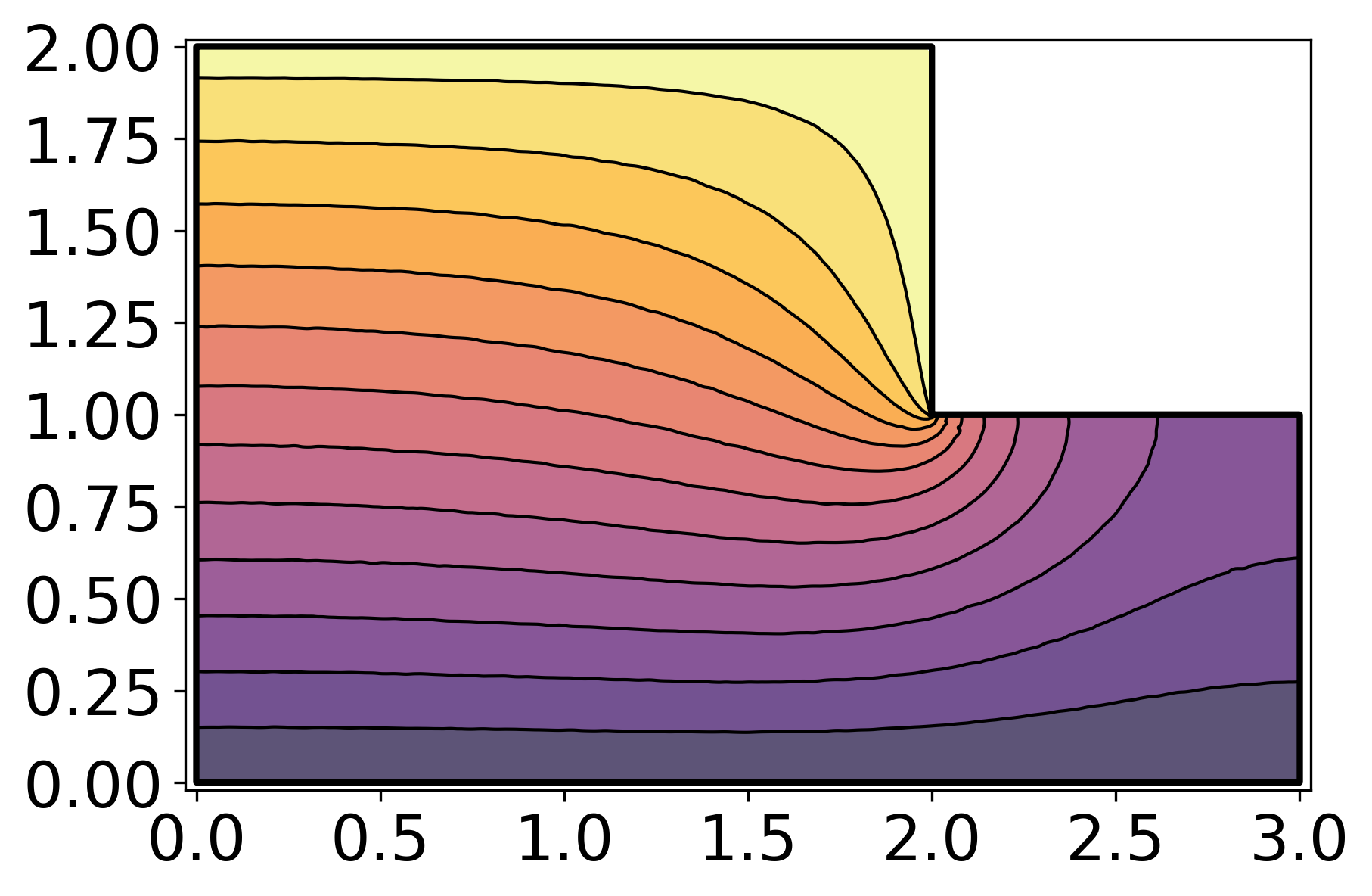}
\includegraphics[width=6.2cm]{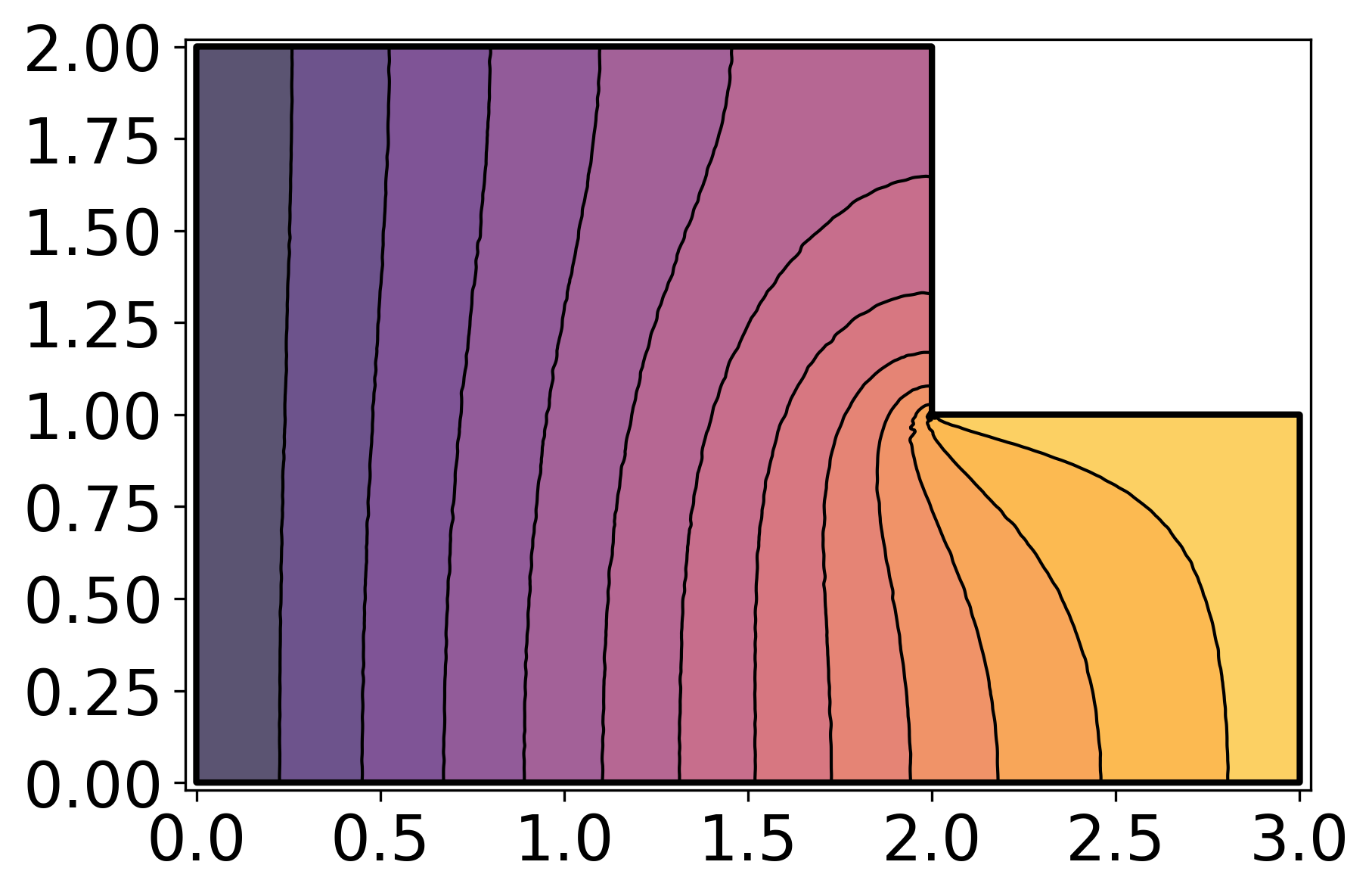}
\caption{Contour lines of the potential functions in the L-shaped region. The vertices of the region $Q$ are \(z_1 = (0, 0)\), \(z_2 = (3, 0)\), \(z_3 = (3, 1)\), \(z_4 = (2, 1)\), \(z_5 = (2, 2)\), and \(z_6 = (0, 2)\). Modulus reference for \(M(Q; z_2, z_4, z_6, z_1) = 1.508154\) and \(M(Q; z_1, z_2, z_4, z_6) = 0.663062\). Computed modulus using the five-point formula is \(1.508108\) and \(0.663082\).}
\label{fig:L_shaped_potential}
\end{figure}

\begin{table}[H]
\centering
\begin{tabular}{|c|c|c|c|}
\hline
Nodes & Reference \cite{HakulaRasilaVuorinen2011} & Computed value & Error \\
\hline
$M(Q; z_2, z_4, z_6, z_1)$ & 1.508154 & 1.508108 & $4.6 \cdot 10^{-5}$ \\
$M(Q; z_1, z_2, z_4, z_6)$ & 0.663062 & 0.663150 & $8.8 \cdot 10^{-5}$ \\
\hline
\end{tabular}
\caption{Moduli of the L-shaped region and its conjugate quadrilateral.}
\label{table:L_shaped_region}
\end{table}

\begin{figure}[H]
\centering
\includegraphics[width=0.55\textwidth]{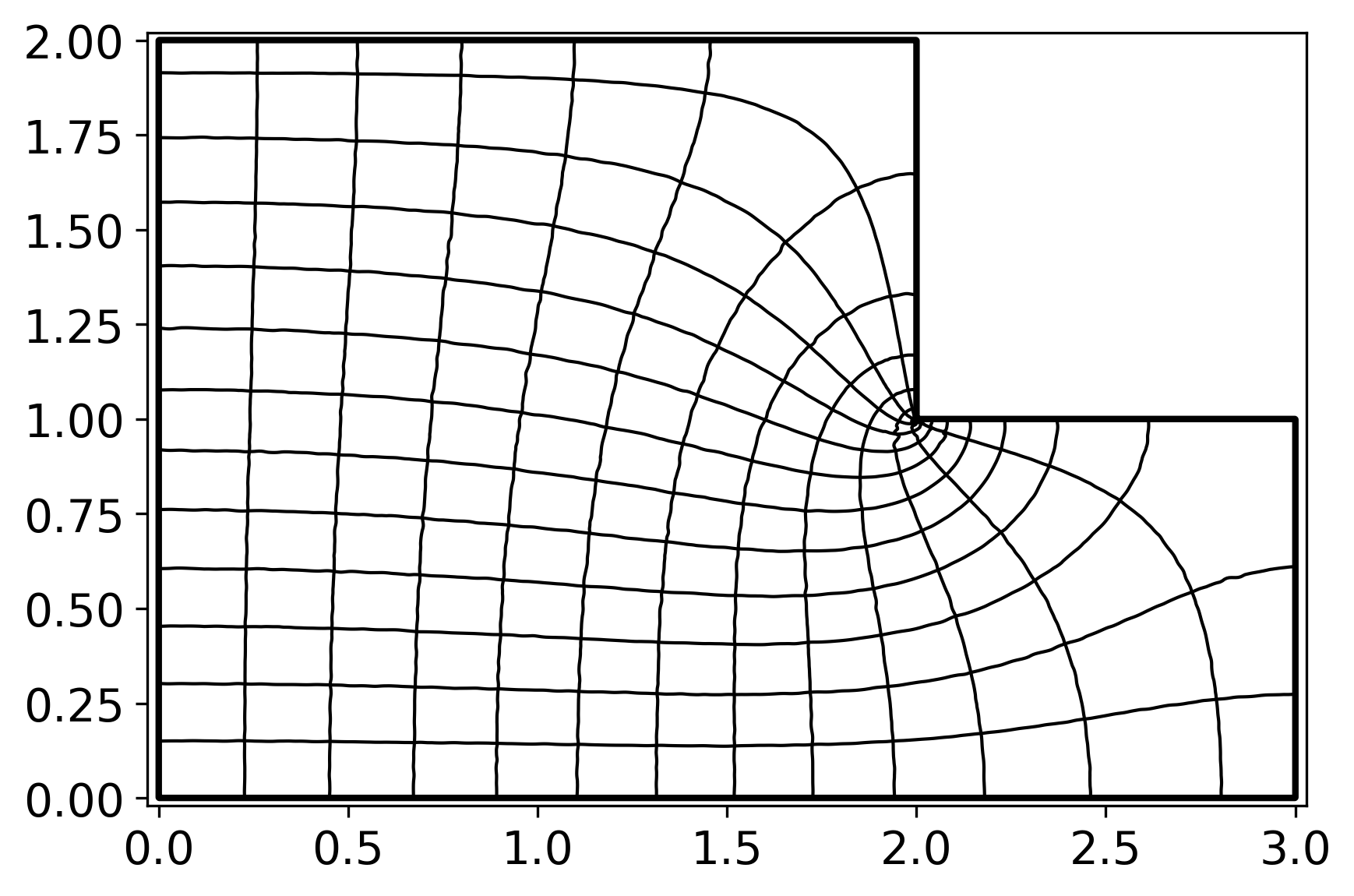}
\caption{The canonical conformal mesh of the L-shaped region.}
\label{fig:L_shaped_mesh}
\end{figure}

\subsection{Circular-arc quadrilaterals}

Given four distinct points \( z_1, z_2, z_3, z_4 \) in the complex plane \( \mathbb{C} \), the cross-ratio is defined by the expression:
\begin{equation}
    \left[z_1, z_2, z_3, z_4\right] = \frac{(z_1 - z_3)(z_2 - z_4)}{(z_1 - z_4)(z_2 - z_3)}.
\end{equation}

An important property of the cross-ratio is its invariance under Möbius transformations. That is, if \( T \) is a Möbius transformation, then the cross-ratio of the transformed points \( T(z_1), T(z_2), T(z_3), T(z_4) \) remains unchanged:

\begin{equation}
    \left[T(z_1), T(z_2), T(z_3), T(z_4)\right] = [z_1, z_2, z_3, z_4].
\end{equation}

\textbf{Type A.} Consider a quadrilateral with sides formed by circular arcs from intersecting orthogonal circles, where angles are \(\pi / 2\). Let \(0<a<b<c<2 \pi\) and select points \(\{1, e^{ia}, e^{ib}, e^{ic}\}\) on the unit circle. The absolute value of the cross-ratio of these four points is

\begin{equation}
\Big| \big[ 1, e^{i a}, e^{i b}, e^{i c} \big] \Big| = \frac{\sin(b/2)\sin((c - a)/2)}{\sin(a/2)\sin((c - b)/2)} = u.
\label{eq:AbsoluteRatioForFourPointsInUnitCircle}
\end{equation}

\begin{figure}[H]
\centering
\includegraphics[width=6cm]{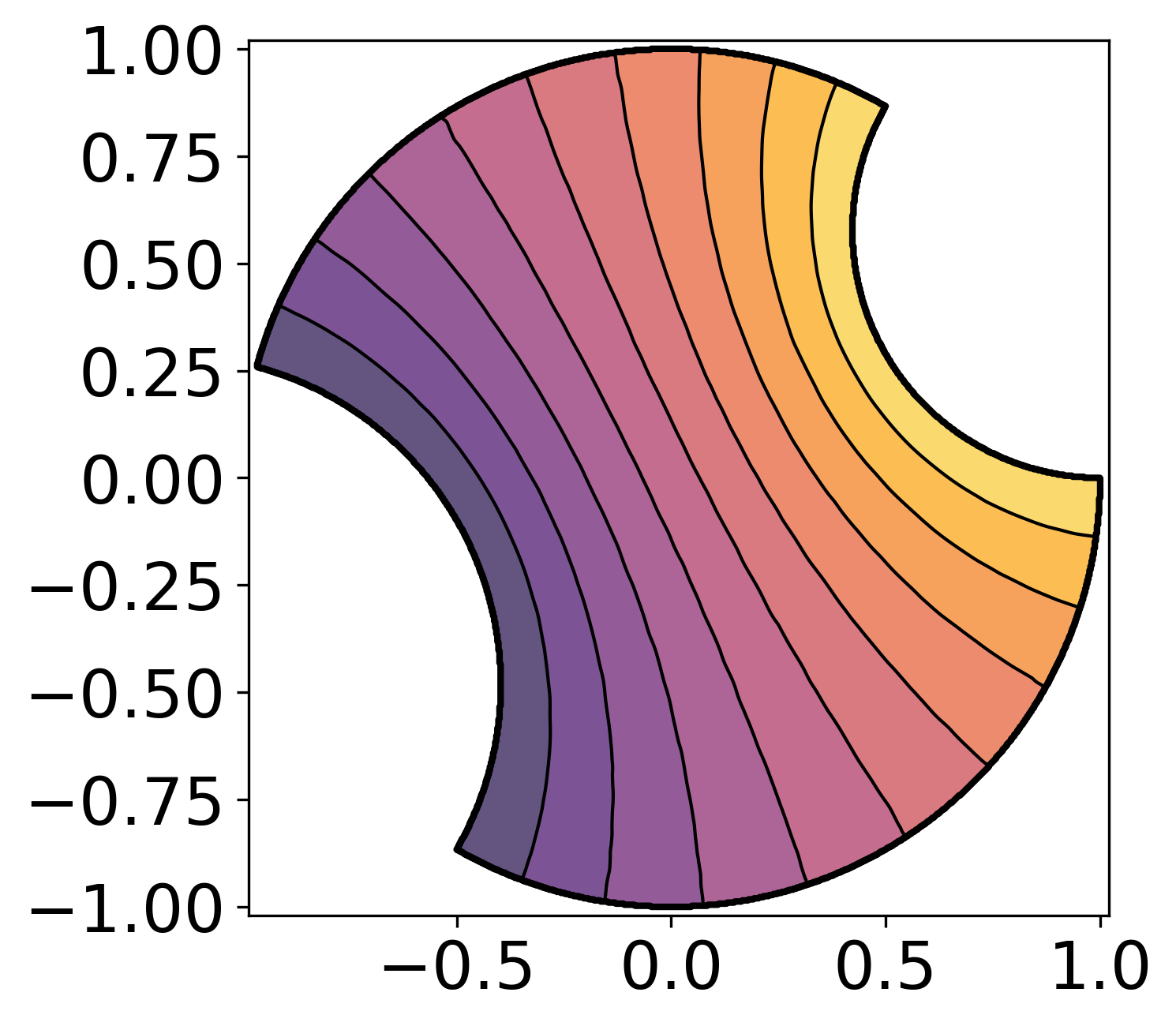}
\includegraphics[width=6cm]{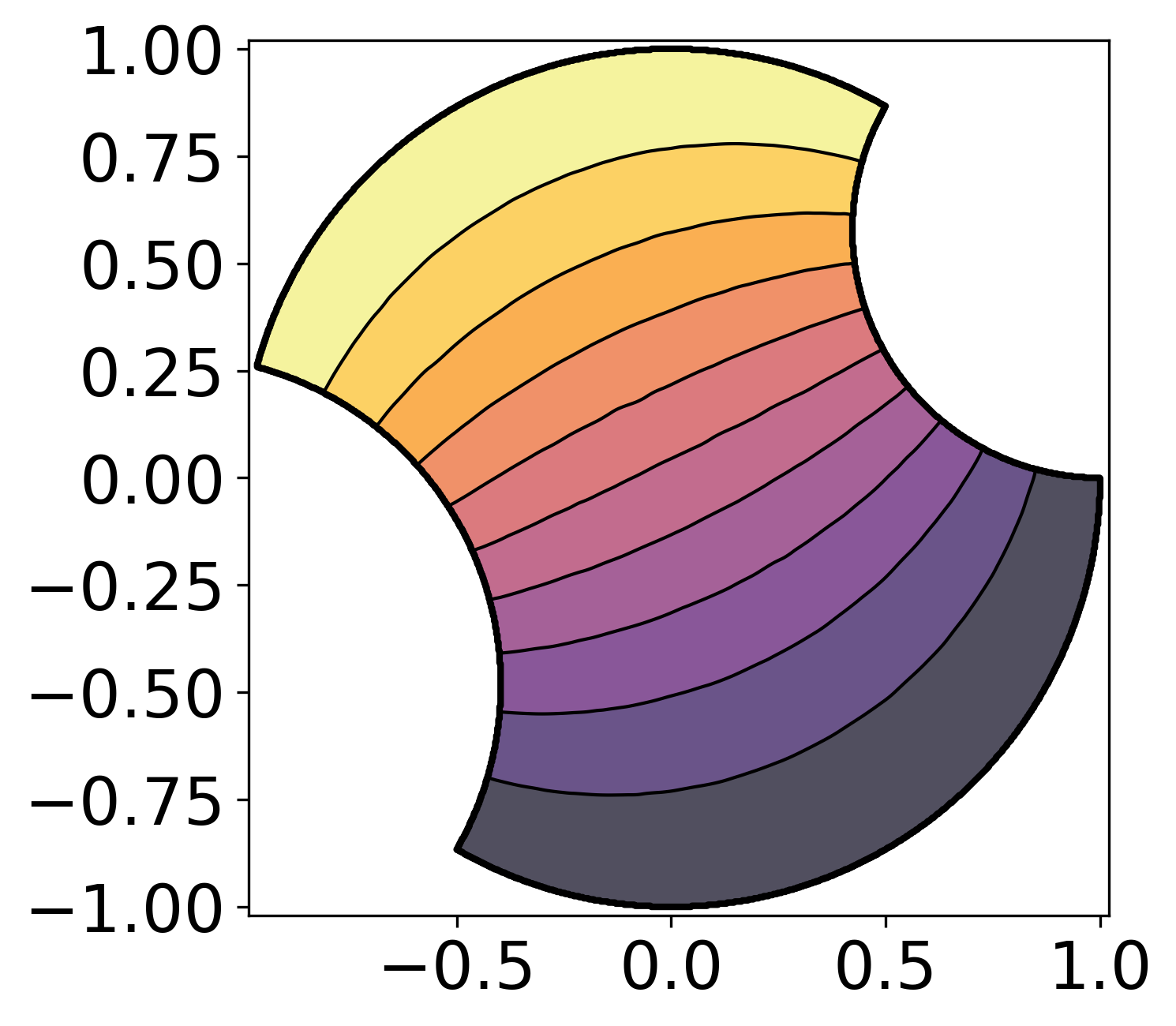}
\caption{The quadrilateral \((Q_A; \pi / 3, 11\pi / 12, 4 \pi / 3, 1)\).}
\label{fig:circular_quadrilaterals}
\end{figure}

Define \(Q_A\) as the domain obtained from the unit disk by removing regions bounded by two orthogonal arcs with endpoints \(\{1, e^{ia}\}\) and \(\{e^{ib}, e^{ic}\}\) respectively, forming the quadrilateral \((Q_A; e^{ia}, e^{ib}, e^{ic}, 1)\).
By applying an appropriate Möbius transformation, we are able to conformally map the quadrilateral \( Q_A \) onto the upper semicircle of the annulus \( \{z \in \mathbb{C} : 1 < |z| < t\} \) in the complex plane. The conformal modulus of the annulus \( \{z \in \mathbb{C} : 1 < |z| < t\} \) is $(2\pi / \log t)$ \cite{Gaier1995}. Therefore the modulus of the quadrilateral is precisely half that of the full annulus,
\begin{equation}
\mathcal{M}(Q_A; e^{i a}, e^{i b}, e^{i c}, 1) = \frac{\pi}{\log t},
\end{equation}

with:
\[
t = 2u - 1 + 2\sqrt{u^2 - u}, \quad t > 1.
\]
The contours of the two potential functions are shown in Figure~\ref{fig:circular_quadrilaterals}.
The computed values of the moduli are in Table~\ref{table:circular_quadrilaterals_moduli}.

\begin{figure}[H]
\centering
\includegraphics[width=0.7\textwidth]{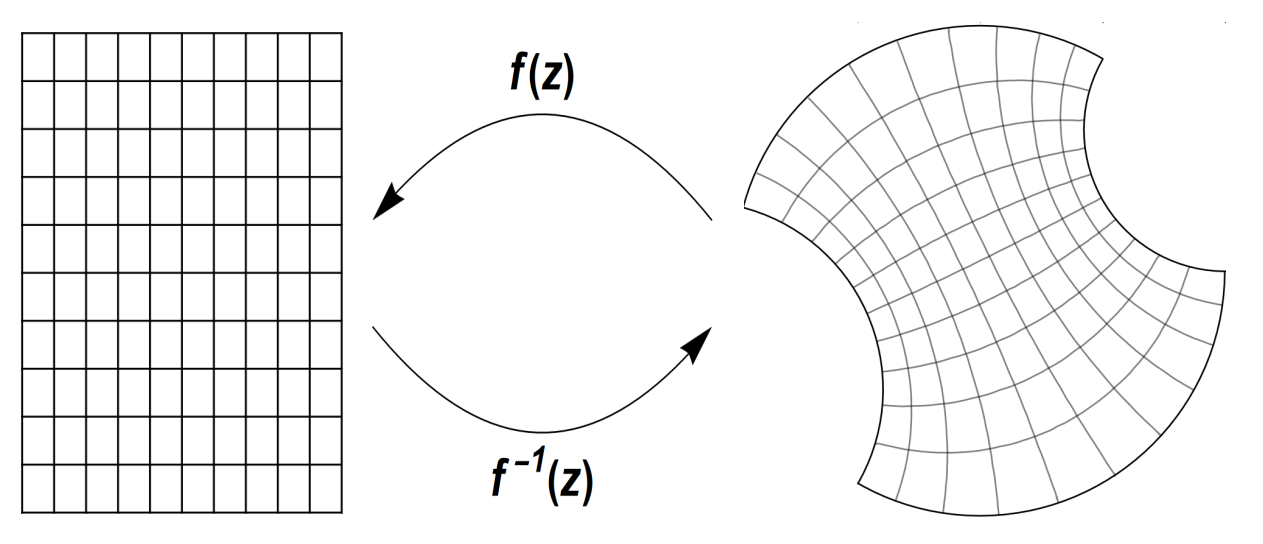}
\caption{The canonical conformal mesh of \((Q_A; \pi / 3, 11\pi / 12, 4 \pi / 3, 1)\).}
\label{fig:circular_quadrilateral_mesh}
\end{figure}

\begin{table}[!htbp]
\centering

\begin{tabular}{|c|c|c|c|}
\hline
Nodes & Reference \cite{HakulaRasilaVuorinen2011} & Computed value & Error \\
\hline
(2, 10, 12) & 0.707150 & 	0.708127 & $9.77 \cdot 10^{-4} $ \\
(2, 10, 14) & 0.807451 &  0.807155 & $2.96 \cdot 10^{-4}$ \\
(4, 12, 18) & 1.038325 &  1.038780 & $4.55 \cdot 10^{-4}$ \\
(6, 16, 24) & 1.170060  &  1.170932 & $8.72 \cdot 10^{-4}$ \\
(8, 22, 32) & 1.313262  &  1.313680 & $4.18 \cdot 10^{-4}$ \\
\hline
\end{tabular}
\caption{Moduli of quadrilaterals \( (Q_{A}; e^{im\pi/24}\), \(e^{in\pi/24}\), \(e^{ir\pi/24}, 1) \) for several integer triples \( (m, n, r) \)}
\label{table:circular_quadrilaterals_moduli}
\end{table}

\textbf{Type B.} 
Consider the sides of the quadrilateral as circular arcs within the unit disk, where each internal angle is \(\pi\). Consequently, the unit disk, along with the boundary points \(\{ e^{ia}, e^{ib}, e^{ic}, 1\}\), defines a quadrilateral, which we denote as \(Q_B\).
Again, by employing a proper Möbius transformation to map the unit disk onto the upper half-plane, we can conveniently express the modulus using the capacity of the Teichmüller ring domain \cite[Section 7]{Anderson1997}. This is formulated as:

\begin{equation}
\mathcal{M}(Q_B; e^{ia}, e^{ib}, e^{ic}, 1) = \frac{1}{2} \tau(u - 1),
\end{equation}

where \(u\) is defined as in equation \eqref{eq:AbsoluteRatioForFourPointsInUnitCircle}, and 

\[
\tau(t) = \frac{\pi}{\mu_{1/2}(1/\sqrt{1 + t})}, \quad t > 0,
\]

\begin{figure}[H]
\centering
\includegraphics[width=6.2cm]{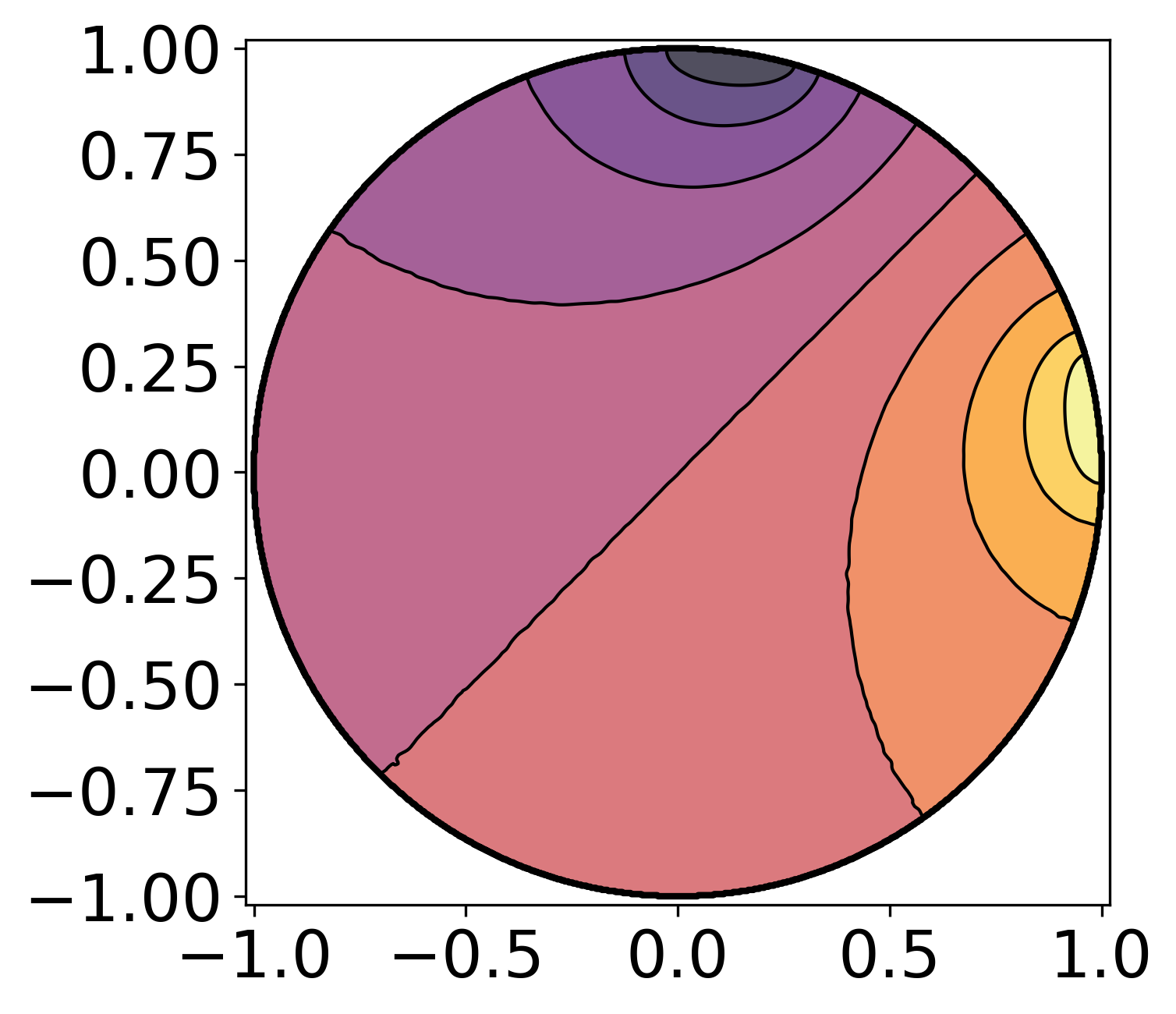}
\includegraphics[width=6.2cm]{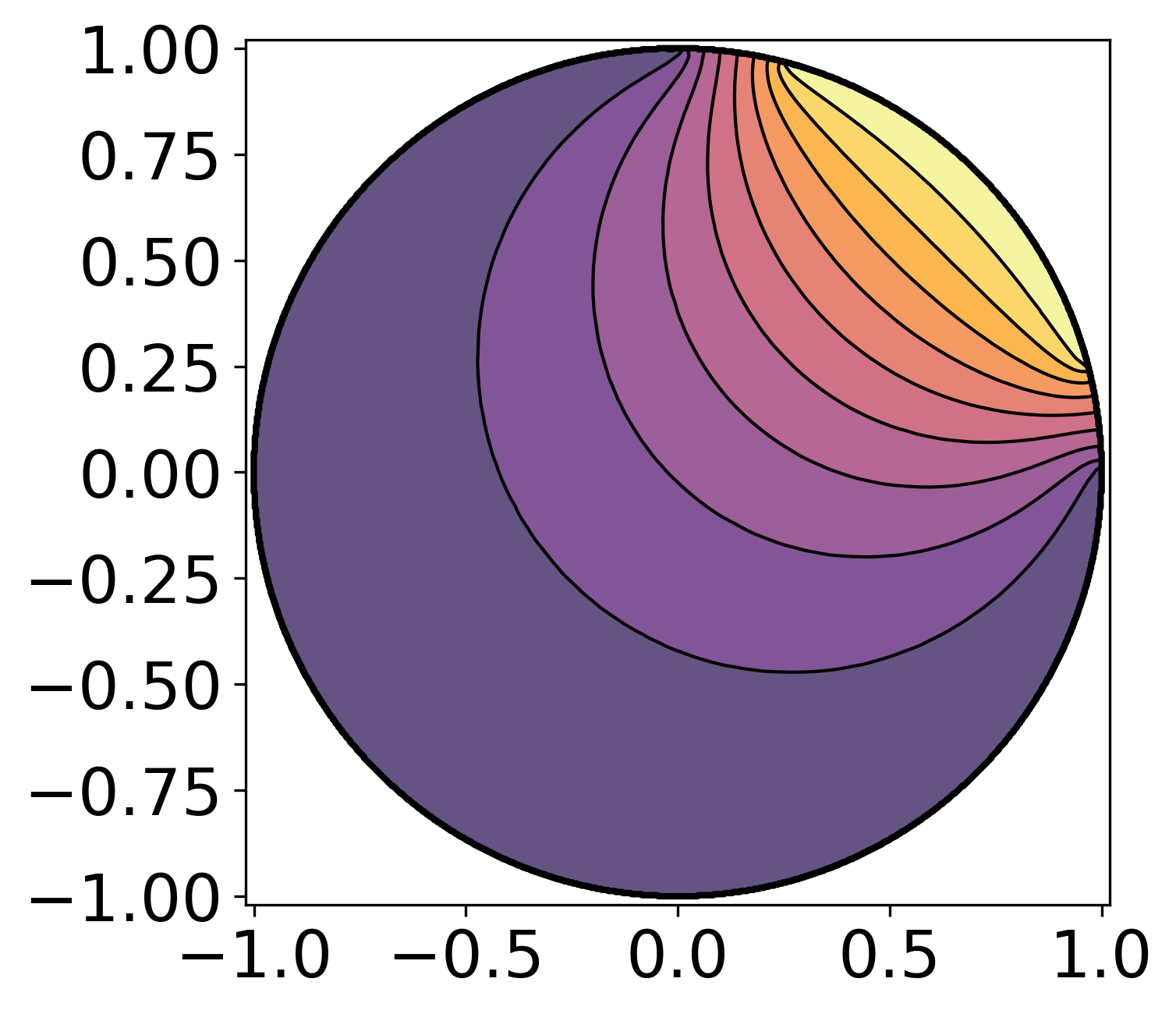}
\caption{The quadrilateral \((Q_B; \pi / 12, 5\pi / 12,  \pi / 2, 1)\).}
\label{fig:circular_quadrilaterals_typeB}
\end{figure}

Here, \(\mu_{1/2}(r)\) is defined in  \cite[Subsection 4.3]{HakulaRasilaVuorinen2011}. The summarized results are presented in Table~\ref{table:circular_quadrilaterals_TypeB_moduli}, The contours of an example are shown in Figure~\ref{fig:circular_quadrilaterals_typeB}.

\begin{figure}[H]
\centering
\includegraphics[width=0.5\textwidth]{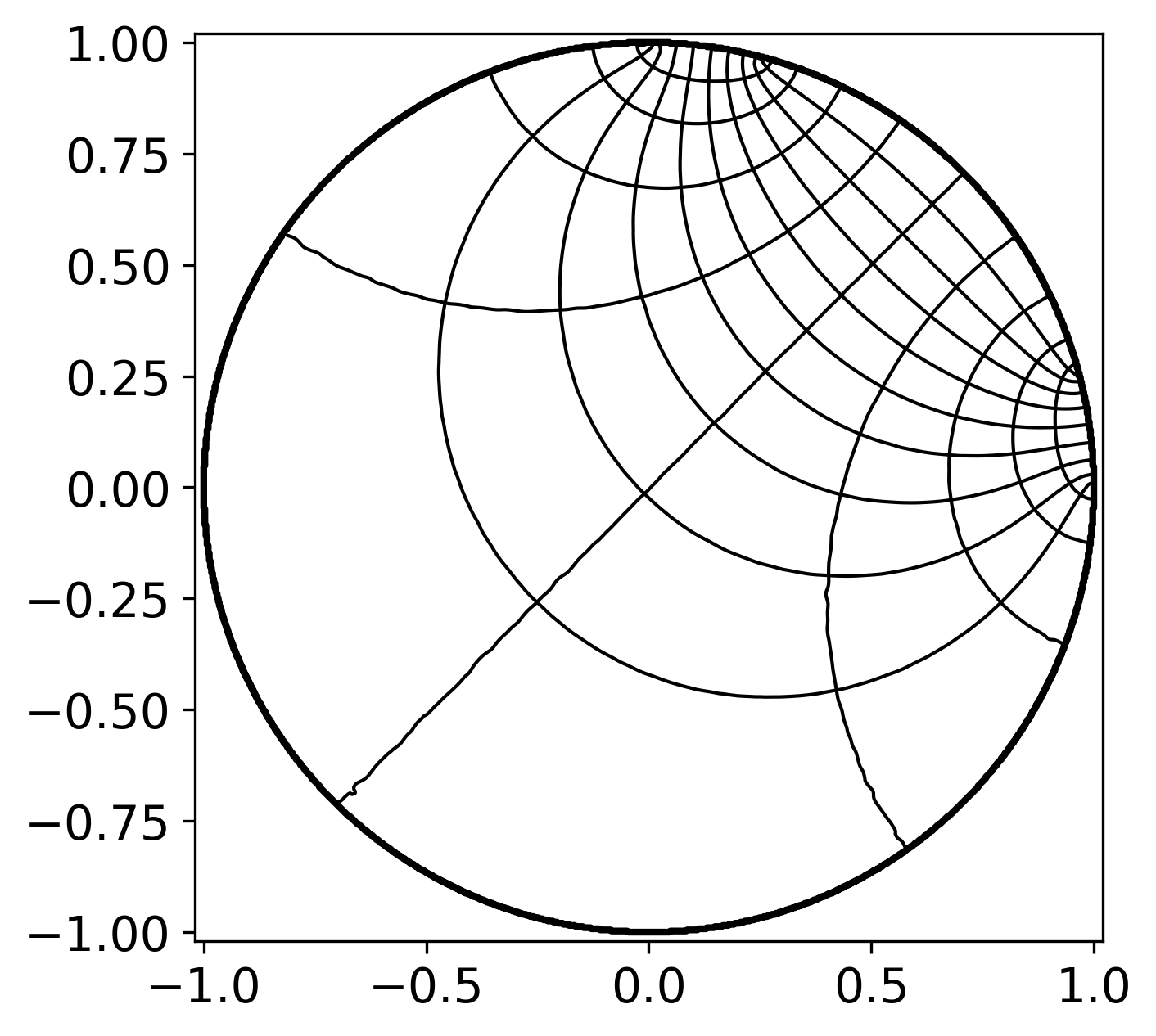}
\caption{The canonical conformal mesh of \((Q_B; \pi / 12, 5\pi / 12,  \pi / 2, 1)\).}
\label{fig:circular_quadrilateral_mesh}
\end{figure}

\begin{table}[!htbp]
\centering

\begin{tabular}{|c|c|c|c|}
\hline
Nodes & Reference \cite{HakulaRasilaVuorinen2011} & Computed value & Error \\
\hline
(2, 10, 12) & 0.538971 & 	0.538686 & $2.85 \cdot 10^{-4} $ \\
(2, 10, 14) & 0.595343 &  0.595158 & $1.85 \cdot 10^{-4}$ \\
(4, 12, 18) & 0.712162 &  0.712533 & $3.71 \cdot 10^{-4}$ \\
(6, 16, 24) & 0.771869  & 0.771132 & $7.37 \cdot 10^{-4}$ \\
(8, 22, 32) & 0.831900  & 0.831458 & $4.42 \cdot 10^{-4}$ \\
\hline
\end{tabular}
\caption{Moduli of quadrilaterals \( (Q_{B}; e^{im\pi/24}\), \(e^{in\pi/24}\), \(e^{ir\pi/24}, 1) \) for several integer triples \( (m, n, r) \)}
\label{table:circular_quadrilaterals_TypeB_moduli}
\end{table}

\vfill
\subsection{Three-dimensional examples}\label{sec-3d}

\begin{figure}
\centering
\includegraphics[width=0.8\textwidth]{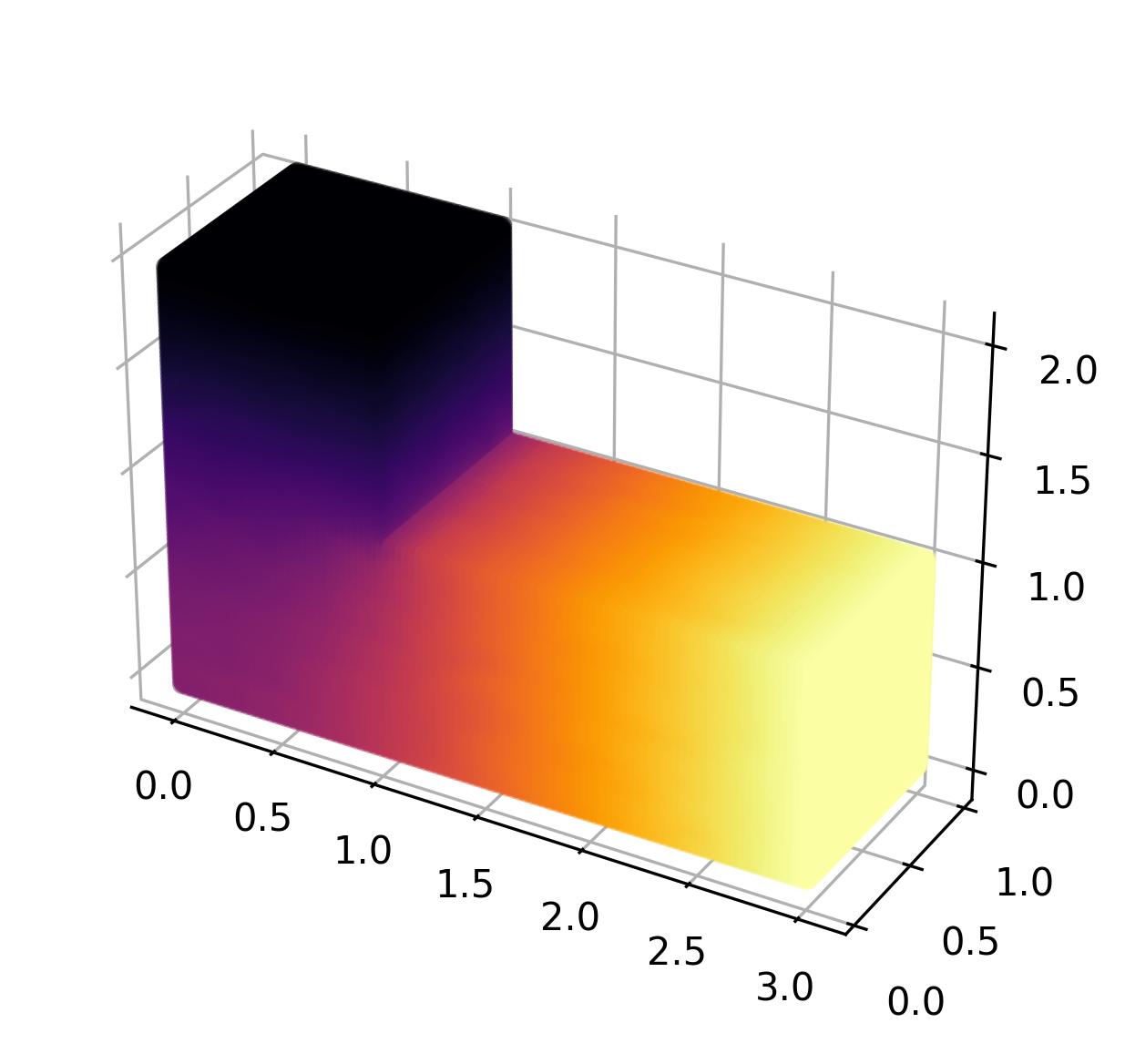}
\caption{Visualization of the solution to Laplace equation in L-shaped polyhedron. The top side and the rightmost side are Dirichlet boundaries with value 1 and 0, respectively. Other sides have zero Neumann boundaries.}
\label{fig:3D_example}
\end{figure}

\begin{figure}
\centering
\includegraphics[width=0.8\textwidth]{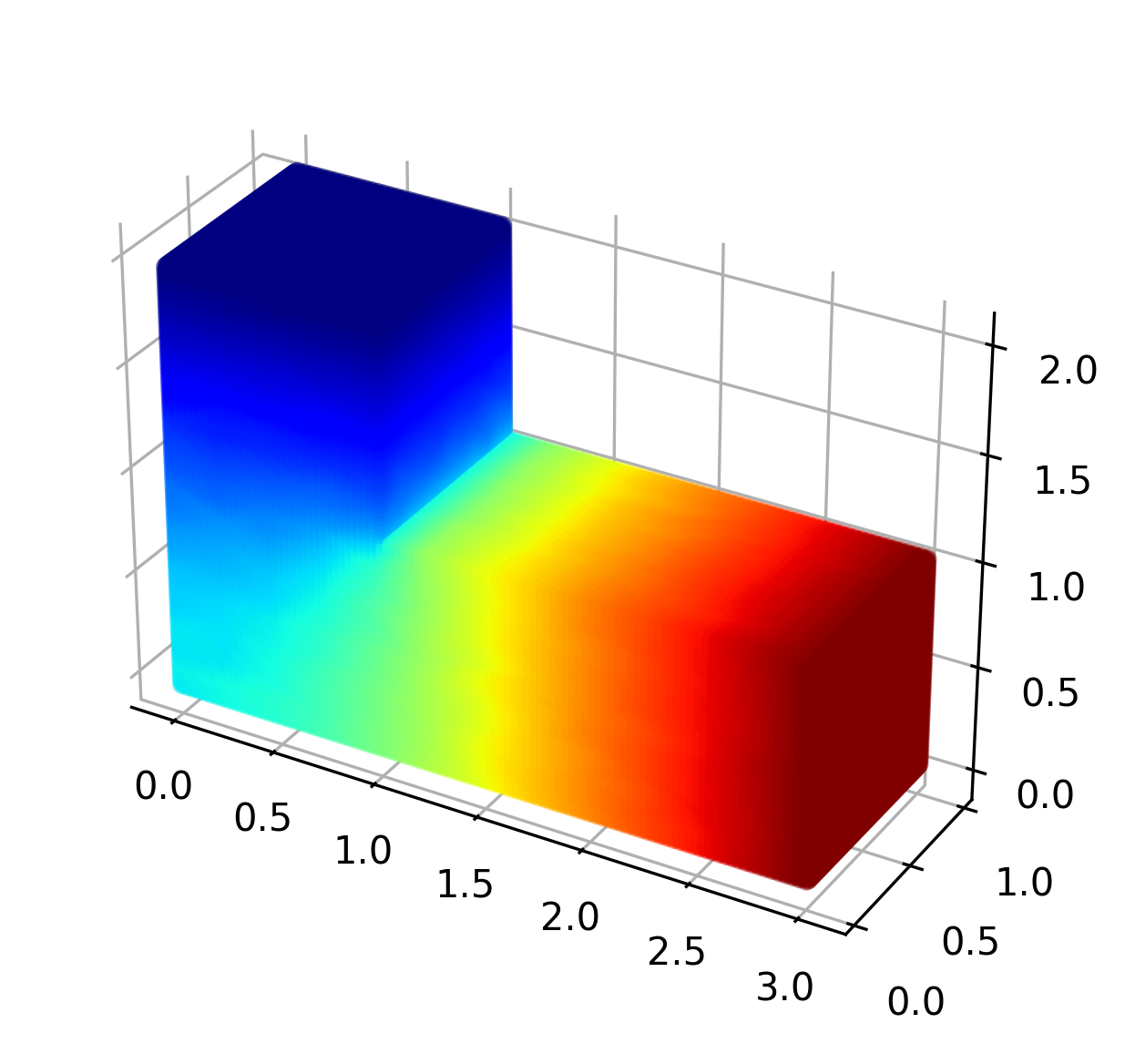}
\caption{(Redrawn in temperature scale) Temperature distribution in a steady heat transfer flow of an L-shaped insulated polyhedron. The top blue side represents the cold sink and the right yellow side represents the hot source.}
\label{fig:3D_example2}
\end{figure}

In higher dimensions, our method has particular appeal due to its relative simplicity when compared to alternative solutions. The reflected WoS algorithm extends naturally to 3D Dirichlet--Neumann domains, modeling insulated heat flow in thermodynamics. See Figures \ref{fig:3D_example} and \ref{fig:3D_example2}.

\section{Discussion}

While several methods for solving the potential functions required by the conjugate function method are known, and the main motivation of this paper was to test feasibility of using a Walk-on-Spheres algorithm for this purpose, this approach does have several distinct advantages. Potential advantages of the approach include the following. 

Firstly, unlike finite element methods, it does not require a mesh or a parameterization of the underlying domain, just a simple representation of the boundary as a polygon or a circular arc polygon is needed. This is a significant advantage in the higher dimensions and in the cases where the domain has a complex geometry \cite{HRV3}. Our algorithm can take a full advantage of parallel computation, because the value of the function is simulated at each point independently of the others. In the planar case, it is suitable for computations on polygonal domains and circular arc polygons, latter of which is not easy to handle in some approaches. The algorithm is also expected to be suitable for unbounded domains (see \cite{HRV2}), although this needs further investigation.

\section*{Code and data availability}

Open source C/Python source code for reproducing results in this paper, along with test results, can be downloaded from the page:\\ \href{https://github.com/arasila/wos_conformal/}{{\tt https://github.com/arasila/wos\_conformal/}}.

\section*{Acknowledgments}
The research was partly supported by NSF of China under the number 11971124, NSF of Guangdong Province under the numbers 2021A1515010326 and 2024A1515010467, and Li Ka Shing Foundation under the number 2024LKSFG06.

Authors wish tho thank Siqi Liu, a Guangdong Technion student, for her help in writing the Python code that was used in plotting the figures in this paper.





\end{document}